\documentclass[10pt,a4paper]{article}

\usepackage{lineno,hyperref}
\modulolinenumbers[5]
\usepackage[utf8]{inputenc}
\usepackage{amsthm,amsmath,enumerate,epic}
\usepackage[mathscr]{euscript}
\usepackage{pstricks,pst-node}
\usepackage{amssymb}
\usepackage{graphicx}
\usepackage{caption}
\usepackage{subcaption}
\usepackage{hyperref}
\usepackage[vcentermath]{youngtab}
\usepackage{parskip}
\theoremstyle{plain}
\newtheorem{thm}[subsubsection]{Theorem}
\newtheorem{lem}[subsubsection]{Lemma}
\newtheorem{prop}[subsubsection]{Proposition}

\newtheorem{ex}[subsubsection]{Example}

\theoremstyle{definition}
\newtheorem{rem}[subsubsection]{Remark}
\newtheorem{defn}[subsubsection]{Definition}

\newcommand\brem{\begin{rem}\begin{sffamily}\begin{upshape}}
\newcommand\erem{\end{upshape}\end{sffamily}\end{rem}}
\newcommand\bdefn{\begin{defn}\begin{rm}}
\newcommand\edefn{\end{rm}\end{defn}}
\newcommand\bex{\begin{ex}\begin{rm}}
\newcommand\eex{\end{rm}\hfill$\Box$\end{ex}}

\newcommand{\bN}{\mathbb{N}}

\newcommand{\ga}{\alpha}
\newcommand{\gb}{\beta}
\newcommand{\gc}{\gamma}

\newcommand{\disjunion}{\,\,\dot{\cup}\,\,}

\newcommand{\gB}{\bar{\gb}\times\gb}
\newcommand{\BBZ}{\gB}

\newcommand{\gtS}{\gtrdot}
\newcommand{\ltS}{\lessdot}

\newcommand{\la}[1]{\stackrel{#1}{\longleftarrow}}

\newcommand\id{I(d)}

\newcommand\modV{(V)}
\newcommand\misd{\mathfrak{M}_d\modV}
\newcommand\miso{\misd}

\newcommand\affinev{\mathbb{A}\modv}

\newcommand\ortho{\ensuremath{\mathfrak{O}}}

\newcommand\dortho{\ensuremath{\mathfrak{A}}}

\newcommand\pos{\mathfrak{N}(v)}
\newcommand\roots{\mathfrak{R}(v)}
\newcommand\posb{\mathfrak{N}(\gb)}
\newcommand\rootsb{\mathfrak{R}(\gb)}

\newcommand\modv{}

\newcommand\st{\,|\,}

\newcommand\posv{\ortho\pos\modv}

\newcommand\rootsv{\ortho\roots\modv}
\newcommand\andposv{\pos\modv}
\newcommand\androotsv{\roots\modv}
\newcommand\dposv{\dortho\pos\modv}

\newcommand\drootsv{\dortho\roots\modv}

\newcommand\orootsb{\ortho\rootsb}

\newcommand\androotsb{\rootsb\modv}

\newcommand\diag{\mathfrak{d}}
\newcommand\diagv{\diag^{v}}

\newcommand\hash{\#}

\newcommand\init{\textup{in}_{\torder}}
\newcommand\torder{\vartriangleright}

\newcommand\card\#
\numberwithin{equation}{subsection}
\numberwithin{figure}{subsection}
\begin{document}
\title{\textbf{Initial ideals of tangent cones to Richardson varieties in the symplectic Grassmannian}}
\date{}
\author {Papi Ray and Shyamashree Upadhyay\\
\\Department of Mathematics\\ Indian Institute of Technology, Guwahati\\Assam-781039, INDIA\\email: popiroy93@iitg.ac.in\\
\\
Department of Mathematics\\ Indian Institute of Technology, Guwahati\\Assam-781039, INDIA\\email: shyamashree@iitg.ac.in}
\maketitle
\begin{abstract}
We give an explicit Gr\"obner
basis for the ideal of the tangent cone at any $T$-fixed point of a Richardson variety in the symplectic Grassmannian, thus generalizing a result of Ghorpade and Raghavan \cite{gr}.
\end{abstract}
\tableofcontents
\section{\textbf{Introduction}}\label{s.introduction}
The study of Schubert varieties has a long and rich history. Richardson varieties are a natural generalization of Schubert varieties. We are interested in Richardson varieties in the symplectic Grassmannian. We consider initial ideals of tangent cones to Richardson varieties in the symplectic Grassmannian. In this paper, we give an explicit Gr\"obner
basis for the ideal of the tangent cone at any $T$-fixed point of a Richardson variety in the symplectic Grassmannian.
\par 
In \cite{Ko-Ra}, Kodiyalam and Raghavan provide (with respect to certain conveniently chosen term orders) an explicit Gr\"obner basis for the ideal of the tangent cone at any $T$-fixed point of a Schubert variety in the ordinary Grassmannian, thereby proving the conjectures of Kreiman and Lakshmibai (made in \cite{Kr-La}). Then in \cite{gr}, Ghorpade and Raghavan do the analogous work for Schubert varieties in the symplectic Grassmannian. And finally in \cite{Ra-Up1,Ra-Up2}, Raghavan and Upadhyay do the analogous work for Schubert varieties in the orthogonal Grassmannian. \par
The above results on Schubert varieties do not admit a straight forward generalization to Richardson varieties. The local properties of Schubert varieties at any $T$-fixed point determine the local properties at all other points, because of the $B$-action; but this does not extend to Richardson varieties, since Richardson varieties only have a $T$-action. However, in \cite{Kr-bkrs}, Kreiman has extended the results of Kodiyalam and Raghavan to Richardson varieties in the ordinary Grassmannian. The analogous work for the orthogonal Grassmannian was done by Upadhyay in \cite{su}. We address the analogous problem for Richardson varieties in the symplectic Grassmannian in this paper.\par
We are motivated by a work of Knutson, Woo and Yong \cite{kwy}, where they give a short proof of the fact that essentially all questions concerning singularities of Richardson varieties reduce to corresponding questions about Schubert varieties. Also, there is a second simplification to the work of \cite{kwy} (One paper to this point is \cite{GK}, due to Graham and Kreiman.). We are also motivated by the method used by Kreiman (in \cite{Kr-bkrs}) to compute an explicit Gr\"obner basis for the ideal of the tangent cone at any $T$-fixed point of a Richardson variety in the ordinary Grassmannian. Our motivation from \cite{kwy} allows us to look at \cite{gr}, where Ghorpade and Raghavan prove that in the case of Schubert varieties in the symplectic Grassmannian, certain objects called ``good admissible pairs'' give rise to a Gr\"obner basis for the ideal of the tangent cone at any $T$-fixed point. In this paper, we have defined ``good admissible pairs'' as a natural extension of the ``good admissible pairs'' of \cite{gr}. Thereafter, we have followed the techniques used by Kreiman in \cite{Kr-bkrs} to obtain an explicit Gr\"obner basis in our case.\par
Sturmfels \cite{stur} and Herzog-Trung \cite{ht} proved results on a class of determinantal varieties which are equivalent to the results of \cite{{Ko-Ra},{kv-thesis},{Kr-bkrs}} for the case of Schubert varieties at the $T$-fixed point $e_{id}$. The key to their proofs was to use a version of the RSK correspondence (see \cite{fulton} for the classical RSK) in order to establish a degree-preserving bijection between a set of monomials defined by an initial ideal and a `standard monomial basis' (see \cite{lr} for a standard monomial basis).\par
In \cite{{Ko-Ra},{kv-thesis}}, an explicit Gr\"obner basis for the ideal of the tangent cone of a Schubert variety in the ordinary Grassmannian at a torus-fixed point were obtained. In \cite{Kr-bkrs}, Kreiman generalizes the results of \cite{{Ko-Ra},{kv-thesis}}   to the case of Richardson varieties. In \cite{Kr-bkrs}, Kreiman gives an explicit Gr\"obner basis for the ideal of the tangent cone at any $T$-fixed point of a Richardson variety in the ordinary Grassmannian, where $T$ denotes a maximal torus in the general linear group. The proof given in \cite{Kr-bkrs} is based on a generalization of the Robinson-Schensted-Knuth (RSK) correspondence, which Kreiman calls the bounded RSK (BRSK). In \cite{papi1}, we had proved that the map $BRSK$ of \cite{Kr-bkrs} and the map $\tilde{\pi}$ of \cite{Ko-Ra} are actually the same maps. In this paper, we use the map $BRSK$ of \cite{Kr-bkrs} to obtain an explicit Gr\"obner basis for the ideal of the tangent cone at any $T$-fixed point of a Richardson variety in the symplectic Grassmannian. The way in which the map $BRSK$ of \cite{Kr-bkrs} has been used here to obtain an explicit Gr\"obner basis has been explained in \S \ref{ss.mainthm+strategy} of this paper.\par
In the study of singularities of Schubert varieties, Woo and Yong investigated Kazhdan-Lusztig ideals \cite{wy1}. These ideals encode coordinates and equations for neighborhoods of type $A$ Schubert varieties at torus fixed points. In \cite{wy2}, Woo and Yong provide a Gr\"obner basis for the Kazhdan-Lusztig ideals. Also in \cite{billey}, the authors discuss three natural generalizations of Richardson varieties which they call projection varieties, intersection varieties, and rank varieties. In \cite{billey}, they study the singularities of each type of generalization.\par
 The organization of the paper is as follows. In \S \ref{ss.In-state-problem}, we define the main objects of interest, namely, the symplectic Grassmannian and Richardson varieties in it. In \S \ref{s.SGRV}, we recall all the things necessary to state the main result of the paper. The main result of the paper comes as Theorem \ref{t.main}, and in \S \ref{ss.mainthm+strategy}, we provide a strategy to prove this theorem. In \S \ref{s.two-sets}, we define the two sets needed to prove the main theorem and then we provide the main proof in \S \ref{s.theproof}. In a forthcoming paper, the results of this paper will be applied to give a combinatorial description of the multiplicity at any torus fixed point of a Richardson variety in the symplectic Grassmannian.
\section{\textbf{Notation and Preliminaries}}\label{s.SGRV}
\subsection{\textbf{\textit{Symplectic Grassmannian and Richardson varieties}}}\label{ss.In-state-problem}
The following definitions and notation are written in the same way as in \cite{gr}. Given any positive integer $n$, we denote by $[n]$ the set $\{1,2,\ldots,n\}$. Given positive integers $r$ and $n$ with $r\leq n$, we denote by $I(r,n)$ the set of all $r$-element subsets of $[n]$. Let $\alpha=(\alpha_1, \ldots , \alpha_r)\in I(r,n)$, where $1\leq\alpha_1< \ldots < \alpha_r\leq n$. If $\beta=(\beta_1, \ldots ,\beta_r)\in I(r,n)$ be such that $1\leq\beta_1< \ldots <\beta_r\leq n$, then we say that $\alpha \leq \beta$ if $\alpha_i \leq \beta_i$ for all $i=1,\ldots,r$. Clearly, $\leq$ defines a partial order on $I(r,n)$.\par
A positive integer $d$ will be kept fixed throughout this paper. For $j \in [2d]$, set $j^* : = 2d+1-j$. Let $I(d)$ denote the set of all $d$-element subsets $v$ of $[2d]$ with the property that exactly one of $j$, $j^*$ belongs to $v$ for every $j \in [d]$. Clearly $I(d)\subseteq I(d,2d)$. In particular, we have the partial order $\leq$ on $I(d)$ induced from $I(d,2d)$.\par
Fix a vector space $V$ of dimension $2d$ over an algebraically closed field of arbitrary characteristic. Fix a \textbf{non-degenerate skew-symmetric bilinear form} $\langle\ ,\ \rangle$ on $V$. Fix a basis $e_1,\ldots,e_{2d}$ of $V$ such that
   $$\langle e_i,e_j\rangle=\left\{\begin{array}{l}
1\ \mbox{if}\ i=j^*\ \mbox{and}\ i<j,\\
-1\ \mbox{if}\ i=j^*\ \mbox{and}\ i>j,\\
0\ \mbox{otherwise}.
\end{array}\right.$$
A linear subspace $W$ of $V$ is said to be \textbf{isotropic} if the form $\langle\ ,\ \rangle$ \textbf{vanishes identically} on it. Let
$${G}_d(V)=\mbox{ the Grassmannian of all } d\mbox{-dimensional subspaces of } V$$
and 
$$\mathfrak{M}_d(V)= \mbox{ the set of all maximal isotropic subspaces of } V.$$
Then $\mathfrak{M}_d(V)$ is a closed subvariety of ${G}_d(V)$ and is called the \textbf{symplectic Grassmannian}.\par
Let $Sp(V)$ denote the group of all linear automorphisms of $V$ that preserve $\langle\ ,\ \rangle$. The elements of $Sp(V)$ that are diagonal with respect to the basis $e_1,\ldots,e_{2d}$ form a maximal torus $T$ of $Sp(V)$. Similarly the elements of $Sp(V)$ that are upper triangular with respect to $e_1,\ldots,e_{2d}$ form a Borel subgroup $B$ of $Sp(V)$ and the elements of $Sp(V)$ that are lower triangular with respect to $e_1,\ldots,e_{2d}$ form a Borel subgroup opposite to $B$ of $Sp(V)$, it is denoted by $B^-$.\par 
The $T$-fixed points of $\mathfrak{M}_d(V)$ are parametrized by $I(d)$ (as explained in \cite[\S 2]{gr} ). The $B$-orbits (as well as $B^-$-orbits) of $\miso$ are naturally indexed by its $T$-fixed points: each
$B$-orbit (as well as $B^-$-orbit) contains one and only one such point. Let $\ga\in\id$ be arbitrary and let $e_\ga$ denote the corresponding $T$-fixed point of $\miso$. The Zariski closure of the $B$ (resp. $B^-$) orbit through $e_\ga$, with canonical reduced
scheme structure, is called a \textbf{Schubert variety} (resp. \textbf{opposite Schubert variety}), and denoted by $X^\ga$ (resp. $X_\ga$).  For $\ga,\gc\in\id$, the scheme-theoretic
intersection $X_\ga^\gc=X_\ga\cap X^\gc$ is called a
\textbf{Richardson variety}. It can be seen easily that the set consisting of all pairs of elements of $\id$ becomes an indexing set for Richardson varieties in $\miso$. It can also be shown that $X_\ga^\gc$ is
nonempty if and only if $\ga\leq \gc$; and that for $\gb\in\id$, $e_\gb\in X_\ga^\gc$
if and only if $\ga\leq\gb\leq\gc$.\par
For the rest of this paper, $\ga,\gb,\gc$ are arbitrarily fixed elements of $\id$ such that $\ga\leq\gb\leq\gc$.
\subsection{\textbf{\textit{$(r,c)$ pairs and monomials}}}\label{ss.basicnot}
For this subsection, let us fix an arbitrary element $v$ of $I(d,2d)$. We will be dealing extensively with ordered pairs $(r,c)$,
$1\leq r,c\leq 2d$,  such that $r$ is not and $c$ is an entry of~$v$.
Let  $$\androotsv \mbox{ denote the set of all such ordered pairs, that is}, \androotsv=\{(r,c)|r\in\{1,\ldots,2d\}\setminus v, c\in v\}.$$ Set $$\andposv:= \{(r,c)\in\androotsv\st r>c\},$$ $$\rootsv:=\{(r,c)\in\androotsv\st r\leq c^*\},$$ $$\posv:=\{(r,c)\in\androotsv\st r>c, r\leq c^*\}=\rootsv\cap\andposv,$$ $$\diagv:=\{(r,c)\in\androotsv\st r=c^*\},$$ $$\drootsv:=\{(r,c)\in\androotsv\st r\geq c^*\},$$ $$\dposv:=\{(r,c)\in\androotsv\st r>c, r\geq c^*\}.$$ We will refer to $\diagv$ as the \textbf{diagonal}. Let $mon\mathfrak{N}(v)$ denote the set of all monomials in $\andposv$.\par
Figure~\ref{figure.1} in Example \ref{ex.eta} below gives a pictorial look of the above sets.
\begin{ex}\label{ex.eta}
{\rm 
 Let $d=7$ and $v=(1,3,4,7,9,10,13)$.

\vspace{ 1 cm}

\begin{figure}[h]

\setlength{\unitlength}{0.3mm}   
\centering

\begin{picture}(60,100)(0,0)
\matrixput(0,0)(15,0){7}(0,15){7}{\circle*{1}}
\multiputlist(-15,0)(0,15){14,12,11,8,6,5,2}
\multiputlist(0,105)(15,0){1,3,4,7,9,10,13}

\linethickness{.2pt}
\dottedline{0.1}(0,90)(0,75)(30,75)(30,45)(45,45)(45,30)(75,30)(75,0)(90,0)
\put(-2,-2){$\bullet$}
\put(13,13){$\bullet$}
\put(28,28){$\bullet$}
\put(43,43){$\bullet$}
\put(58,58){$\bullet$}
\put(73,73){$\bullet$}
\put(88,88){$\bullet$}
\put(75,5){$\longleftarrow$}
\put(95,5){ boundary of $\mathfrak{N}(v)$}
\end{picture}
\caption{Illustration of the grid representing $\androotsv$}
\label{figure.1}
\end{figure}
The points (including the dark circles) of the above grid represent the set $\mathfrak{R}(v)$ for $v=(1,3,4,7,9,10,13)$. The path sketched on the grid by some piecewise line segments denote the boundary of $\mathfrak{N}(v)$. The points on this grid which lie on the boundary of $\mathfrak{N}(v)$ or to the left of it belong to the set $\mathfrak{N}(v)$. The dark circles denote the diagonal elements. Points above and on the diagonal belong to the set $\mathfrak{OR}(v)$. Again, points which are on and towards the left of the boundary of $\mathfrak{N}(v)$, and which also lie on and above the diagonal, are the points of $\mathfrak{ON}(v)$.}
\end{ex}
We will be considering {\em monomials} in some of these sets. A \textbf{monomial}, as usual, is a subset with each member being allowed a multiplicity (the multiplicity taking values in the non-negative integers). The \textbf{degree} of a monomial also has the usual meaning: consider the underlying set of the monomial and look at the multiplicity with which each element of this underlying set appears in the monomial, the degree of the monomial is the sum of these multiplicities.\par
For a monomial $\mathfrak{S}\in\androotsv$, let $\mathfrak{S}^{\hash}$ denote the set $\{(c^*,r^*)|(r,c)\in\mathfrak{S}\}$.
\begin{ex}
{\rm
Let $d$, $v$ be as in Example~\ref{ex.eta} above. Let $\mathfrak{S}=\{(2,1),(6,4)^2,(5,10)\}$. Then $\mathfrak{S}$ is a monomial in $\mathfrak{R}(v)$. The {underlying set} of the monomial $\mathfrak{S}$ is $\{(2,1),(6,4),(5,10)\}$. The {degree} of $\mathfrak{S}$ is $4$. The {multiplicity} of $(2,1), \ (6,4)$, and $(5,10)$ are respectively $1,\ 2$, and $1$. Also for the monomial $\mathfrak{S}$, $\mathfrak{S}^\# = \{(14,13),(11,9)^2,(5,10)\}$.}
\end{ex}

Let $S$ be any set.  A \textbf{multiset} $E$ on $S$ is defined to be a function $E:S\rightarrow \{0,1,2,\ldots,\}$.  One should think of $E$ as consisting of the set $S$ of elements, but with each $s\in S$ occurring $E(s)$ times. Note that a set is a special type of multiset in which each element occurs exactly once. We call $E(s)$ the \textbf{multiplicity} of $s$
in $E$. Define the multiset $E\disjunion F$ as follows:
\begin{align*}
(E\disjunion F)(s)&=E(s)+F(s),\ s\in S
\end{align*}
Let $\bN$ denote the set of all positive integers. Let $A=\{a_1,a_2,\ldots\}$ and $B=\{b_1,b_2,\ldots\}$ be two multisets on $\bN$ of the same degree, with $a_i\leq a_{i+1}$, $b_i\leq b_{i+1}$, for all $i$. We say that $A$ is less than or equal to $B$ in the \textbf{termwise order} if $a_i\leq b_i$ for all $i$. We denote this by $A\leq B$. We say that $A$ is less than $B$ in the \textbf{strict termwise order} if $a_i< b_i$ for all $i$.  We denote this by $A\ltS B$.\par If $A$, $B$, $C$, and $D$ are multisets on $\bN$ such that
$|A\disjunion D|=|B\disjunion C|$, then we write
\begin{equation}\label{e.m.set_subtraction}
A-C\leq B-D \hbox{ to indicate that }A\disjunion D\leq B\disjunion
C.
\end{equation}
Let  $U=\{(e_1,f_1),(e_2,f_2),\ldots\}$ be a multiset on $\bN^2$.
Define $U_{(1)}$ and $U_{(2)}$ to be the multisets
$\{e_1,e_2,\ldots\}$ and $\{f_1,f_2,\ldots\}$ respectively on
$\bN$. Define the \textbf{non vanishing}, \textbf{negative}, and
\textbf{positive parts} of $U$ to be the following multisets:
\begin{align*}
U^{\neq 0}&=\{(e_i,f_i)\in U\mid e_i-f_i\neq 0\},\\
U^-&=\{(e_i,f_i)\in U\mid e_i-f_i<0\},\\
U^+&=\{(e_i,f_i)\in U\mid e_i-f_i>0\}.
\end{align*}
We say that $U$ is \textbf{non vanishing} if $U\subset (\mathbb{N}^2)^{\neq
0}$, \textbf{negative} if $U\subset (\mathbb{N}^2)^-$, and
\textbf{positive} if $U\subset (\mathbb{N}^2)^+$.  Impose the following transitive relation on multisets on $\bN^2$:
\begin{equation}
U\leq V\iff U_{(1)}-U_{(2)}\leq V_{(1)}-V_{(2)}.
\end{equation}
A \textbf{chain} in $\mathbb{N}^2$ is a subset $C=\{(e_1,f_1),\ldots,(e_m,f_m)\}$ of $\mathbb{N}^2$ such that $e_1<\cdots<e_m$ and $f_1>\cdots>f_m$. Let $T$ and $W$ be negative and positive subsets of $\mathbb{N}^2$ respectively. A non vanishing multiset $U$ on $\mathbb{N}^2$ is said to be \textbf{bounded by} $T,W$ if for every chain $C$ which is contained in the underlying set of $U$, we have
$$T\leq C\leq W.$$
Let $\iota$ be the map on multisets on $\mathbb{N}^2$ defined by, 
$$\iota(\{(e_1,f_1),(e_2,f_2),\ldots\}):=\{(f_1,e_1),(f_2,e_2),\ldots\}.$$
Then $\iota$ is an involution, and it maps negative multisets on $\mathbb{N}^2$ to positive ones and visa-versa.\par
\subsection{\textbf{\textit{Notched and semistandard Young tableaux}}}\label{ss.BRSK}
In this subsection we are recalling the following things from \cite{Kr-bkrs}.\par
A \textbf{Young diagram} (resp. \textbf{notched diagram}) is a collection of boxes arranged into a left and top justified array (resp. into left justified rows). The \textbf{empty Young diagram} is the Young diagram with no boxes. A notched diagram may contain rows with no boxes; however, a Young diagram may not, unless it is the empty Young diagram. A \textbf{Young tableau} (resp. \textbf{notched tableau}) is a filling of the boxes of a Young diagram (resp. notched diagram) with positive integers. The \textbf{empty Young tableau} is the Young tableau with no boxes. Let $P$ be either a notched tableau or a Young tableau. We say that $P$ is \textbf{row strict} if the entries of any row of $P$ strictly increase as one moves to the right. If $P$ is a Young tableau, then we say that $P$ is \textbf{semistandard} if it is row strict and the entries of any column weakly increase as one moves down.\par
Example \ref{eg.rowstrict-ss} below illustrates a row strict notched tableau and a semistandard Young tableau.
\begin{ex}\label{eg.rowstrict-ss}
{\rm
A row strict notched tableau $P$ and a semistandard Young tableau $R$.
$$ \quad P=\young(1234,24,567,4578) \quad \mbox{ and } \quad  R=\young(1234,2456,35,6).$$}
\end{ex}
Let $P$ be a row strict notched tableau and $b$ be a positive integer. Since $P$ is row strict, its entries which are greater than or equal to $b$ are right justified in each row. If we remove these entries (which are greater than or equal to $b$) from $P$, then we are left with a row strict notched tableau, which we denote by $P^{<b}$. We say that $P$ is \textbf{semistandard on} $b$ if $P^{<b}$ is a semistandard Young tableau.\par
Example \ref{eg.P<b} below illustrates $P^{<b}$, for  fixed values of $P$ and $b$.
\begin{ex}\label{eg.P<b}
{\rm
For the row strict notched tableau $P$ in Example~\ref{eg.rowstrict-ss} and $b=4$, we have
$$P^{<b}=\young(123,2).$$
However, for the same $P$, if we take $b=6$, then
$$P^{<b}=\young(1234,24,5,45).$$
Hence $P$ is semistandard on $4$, but not on $6$.}
\end{ex}
A \textbf{notched bitableau} is a pair $(P, Q)$ of notched
tableaux of the same shape (i.e., the same number of rows and the
same number of boxes in each row).  The \textbf{degree} of $(P,Q)$
is the number of boxes in $P$ (or $Q$). A notched bitableau
$(P,Q)$ is said to be \textbf{row strict} if both $P$ and $Q$ are
row strict. A row strict notched bitableau $(P,Q)$ is said to be
\textbf{semistandard} if
\begin{equation}\label{e.s.ss_tabl_1}
P_1-Q_1\leq\cdots\leq P_r-Q_r,
\end{equation}
where $r$ is the total number of rows in $P$ (or $Q$) and for each $i\in\{1,\ldots,r\}$, $P_i$ (resp. $Q_i$) denotes the $i$-th row (from the top) of $P$ (resp. $Q$). A row strict notched bitableau $(P,Q)$ is said to be
\textbf{negative} if $P_i\ltS Q_i$, $i=1,\ldots,r$,
\textbf{positive} if $P_i\gtS Q_i$, $i=1,\ldots,r$, and
\textbf{non vanishing} if either
\begin{equation}\label{e.s.ss_tabl_2}
P_i\ltS Q_i\ \ \text{   or   }\ \ P_i\gtS Q_i,
\end{equation}
for each $i=1,\ldots,r$.\par
Example \ref{eg.ss-notched-bitableau} below gives an illustration of a row strict semistandard non vanishing bitableau.
\begin{ex}\label{eg.ss-notched-bitableau}
{\rm
Consider the notched bitableau
$$(P,Q)=\left(\ \young(123,4567)\ \ \ ,
\ \ \ \young(789,2345)\right).$$ We have that
\begin{enumerate}
\item $(P,Q)$ is row strict.

\item $P_1\disjunion Q_2=\{1,2,2,3,3,4,5\}\leq
\{4,5,6,7,7,8,9\}=P_2\disjunion Q_1$. Therefore, $P_1-Q_1\leq
P_2-Q_2$. Thus $(P,Q)$ is semistandard.

\item $P_1\ltS Q_1$ and $P_2\gtS Q_2$. Thus $(P,Q)$ is non vanishing.
\end{enumerate}}
\end{ex}
Let $(P,Q)$ be a semistandard notched bitableau. If for subsets $T$ and $W$ of $\mathbb{N}^2$,
$$T_{(1)}-T_{(2)}\leq P_1-Q_1\ \mbox{and}\ P_r-Q_r\leq W_{(1)}-W_{(2)},$$
then we say that $(P,Q)$ is \textbf{bounded by} $T,W$.\par 
If $(P,Q)$ is a non vanishing semistandard notched bitableau, then we define $\iota(P,Q)$ to be the notched bitableau obtained by reversing the order of the rows of $(Q,P)$.\par
\subsection{\textbf{\textit{Schensted insertion and bounded insertion}}}\label{sb insertion}
Let us now recall the \textbf{ordinary Schensted insertion} process from \cite[\S 3]{Kr-bkrs}. It is an algorithm which takes as input a semistandard Young tableau $P$, a positive integer $a$, and produces as output a new semistandard Young tableau with the same shape as $P$ plus one extra box, and with the same entries as $P$ (possibly in different locations) plus one additional entry, namely $a$. To begin, insert $a$ into the first row of $P$, as follows. If $a$ is strictly bigger than all entries in the first row of $P$, then place $a$ in a new box on the right end of the first row, and the insertion process terminates. Otherwise, find the smallest entry of the first row of $P$ which is greater than or equal to $a$, and replace that number with $a$. We say that the number which was replaced was ``bumped" from the first row. Insert the bumped number into the second row in precisely the same way as $a$ was inserted into the first row. This process continues down the rows until, at some point, a number is placed in a new box on the right end of some row, at which point the process terminates.\par
We next describe the \textbf{bounded insertion algorithm}, which
takes as input a positive integer $b$, a notched tableau $P$ which
is semistandard on $b$, and a positive integer $a<b$, and produces
as output a notched tableau which is semistandard on $b$, which we denote by $P\la{b}a$.\par
\noindent{\textbf{Bounded Insertion:}} 
\begin{quote}
{\it
\begin{itemize}
\item[\textbf{Step 1.}] Remove all entries of $P$ which are
greater than or equal to $b$ from $P$, resulting in the
semistandard Young tableau $P^{<b}$.
\item[\textbf{Step 2.}] Insert $a$ into $P^{<b}$ using the
ordinary Schensted insertion process (as described above).
\item[\textbf{Step 3.}] Place the entries of $P$ which were
removed when forming $P^{<b}$ in Step 1 back into the Young
tableau resulting from Step 2, in the same rows from which they
were removed.
\end{itemize}
}
\end{quote}
 Example \ref{eg.bounded-insertion} below gives an illustration of the bounded insertion algorithm.
\begin{ex}\label{eg.bounded-insertion}
{\rm
Let $a=3$ and $b=4$. We compute $P\la{b}a$, where 
$$P=\young(1234,24,567,4578).$$ 
Observe that in Step 1, we obtain $P^{<4}$, where $$P^{<4}=\young(123,2).$$
In Step 2, we insert $a=3$ into $P^{<4}$ using the ordinary Schensted insertion process, to get $$\young(123,23).$$ And finally in Step 3, we obtain $$P\la{4}3=\young(1234,234,567,4578).$$\\}
\end{ex}

\subsection{\textbf{\textit{The bounded RSK correspondence}}}\label{brsk correspon}
We next define the \textbf{bounded RSK correspondence}, $BRSK$, a
function which maps negative multisets on $\mathbb{N}^2$ to negative
semistandard notched bitableaux. Let
$$U=\{(a_1,b_1),\ldots,(a_t,b_t)\}$$ be a negative multiset on
$\mathbb{N}^2$, whose entries we assume are listed in
\textbf{lexicographic order}: (i) $b_1\geq\cdots\geq b_t$ and
(ii) if for any $i\in\{1,\ldots,t-1\}$, $b_i=b_{i+1}$, then
$a_i\geq a_{i+1}$. We inductively form a sequence of notched
bitableaux $(P^{(0)},Q^{(0)})$, $(P^{(1)},Q^{(1)}),$
$\ldots,(P^{(t)},Q^{(t)})$, such that $P^{(i)}$ is semistandard on
$b_i$, $i=1,\ldots,t$, as follows:
\begin{quote}
Let $(P^{(0)},Q^{(0)})=(\emptyset,\emptyset)$ and let $b_0=b_1$. Assume
inductively that we have formed $(P^{(i)}, Q^{(i)})$, such that
$P^{(i)}$ is semistandard on $b_i$, and thus on $b_{i+1}$, since
$b_{i+1}\leq b_i$. Define $P^{(i+1)}=P^{(i)}\la{b_{i+1}}a_{i+1}$.
Since bounded insertion preserves semistandardness on $b_{i+1}$,
$P^{(i+1)}$ is also semistandard on $b_{i+1}$. Let $j$ be the row
number of the new box of this bounded insertion. Define
$Q^{(i+1)}$ to be the notched tableau obtained by placing
$b_{i+1}$ on the \it{left} end of row $j$ of $Q^{(i)}$ (and
shifting all other entries of $Q^{(i)}$ to the \it{right} one box).
Clearly $P^{(i+1)}$ and $Q^{(i+1)}$ have the same shape.
\end{quote}
Then $BRSK(U)$ is defined to be $(P^{(t)},Q^{(t)})$. In the process above, we write $$(P^{(i+1)},Q^{(i+1)})=(P^{(i)},Q^{(i)})\la{b_{i+1}}a_{i+1}.$$
In terms of this notation,
$$BRSK(U)=((\emptyset,\emptyset)\la{b_1}a_1)\cdots\la{b_t}a_t.$$
If $U$ is a positive multiset on $\mathbb{N}^2$, then we define $BRSK(U)$ to be $\iota(BRSK(\iota(U)))$.\par
Example \ref{eg.brsk} below gives an illustration of the map $BRSK$.
\begin{ex}\label{eg.brsk}
{\rm
\newcommand{\J}{11}
\newcommand{\B}{13}
Let $U=\{(2,1),(5,3),(6,4),(6,9),(8,13),(11,13)\}$ be a multiset on $\mathbb{N}^2$. Now 
$$\{(2,1),(5,3),(6,4)\} \subset (\mathbb{N}^2)^+$$ and 
$$\{(6,9),(8,13),(11,13)\} \subset (\mathbb{N}^2)^-.$$
Let $U^+ = \{(2,1),(5,3),(6,4)\}$ and $U^- = \{(6,9),(8,13),(11,13)\}$. So $U=U^+ \cup U^-$. Now after arranging $\iota(U^+)$ in lexicographic order, we have $$\iota(U^+)=\{(4,6),(3,5),(1,2)\}.$$
Let us first apply the map $BRSK$ on $\iota(U^+)$. Then\\
\[
\begin{array}{l@{\hspace{.9cm}}l}
P^{(0)}=\emptyset
&Q^{(0)}=\emptyset\\ \\
P^{(1)}=\emptyset\la{6}4=\young(4)
&Q^{(1)}=\young(6)\\ \\
P^{(2)}=\young(4)\la{5}3=\young(3,4)
&Q^{(2)}=\young(6,5)\\ \\
P^{(3)}=\young(3,4)\la{2}1=\young(13,4)
&Q^{(3)}=\young(26,5)\\ \\
\end{array}
\]
\vspace{.5em}
Therefore\\ $BRSK(U^+)=\iota(BRSK(\iota(U^+)))=\left(\, \young(5,26)\ ,\
\young(4,13)\right)$.\\ \\
After arranging $U^-$ in lexicographic order, we have $U^-= \{(11,13),(8,13),(6,9)\}$. Let us now apply the map $BRSK$ on $U^-$. Then\\
\[
\begin{array}{l@{\hspace{.9cm}}l}
P^{(0)}=\emptyset   &
Q^{(0)}=\emptyset\\  \\
P^{(1)}=\emptyset\la{\B}\J=\young(\J) &Q^{(1)}=\young(\B)\\ \\
P^{(2)}=\young(\J)\la{\B}8=\young(8,\J) &Q^{(2)}=\young(\B,\B)\\ \\
P^{(3)}=\young(8,\J)\la{9}6=\young(6,8\J) &Q^{(3)}=\young(\B,9\B)\\ \\
\end{array}
\]
\vspace{.5em}
Therefore $$BRSK(U^-) = \left(\, \young(6,8\J)\ ,\
\young(\B,9\B)\right)
\mbox{ and } BRSK(U) = \left(\, \young(6,8\J,5,26)\ ,\
\young(\B,9\B,4,13)\right).$$}
\end{ex}
\subsection{\textbf{\textit{The Kodiyalam-Raghavan maps}}}\label{ss.pi}
In this subsection, we will recall the maps $\pi$ and $\tilde{\pi}$  of Kodiyalam-Raghavan \cite{Ko-Ra}. \par 
Fix an element $v$ in $I(d,2d)$. The map $\pi$ is from $mon\mathfrak{N}(v)$ to $I(d,2d)\times mon\mathfrak{N}(v)$, and the map $\tilde{\pi}$ is from $mon\mathfrak{N}(v)$ to $\widetilde{SM^{v,v}}$.\par
Given any ${\beta}_1=(r_1,c_1), \  {\beta}_2=(r_2,c_2)$ in $\mathfrak{N}(v)$, we say that ${\beta}_1=(r_1,c_1) > {\beta}_2=(r_2,c_2)$ if $r_1 >r_2$ and $c_1<c_2$. A sequence ${\beta}_1>\ldots>{\beta}_t$ of elements of $\mathfrak{N}(v)$ is called a $v$-\textbf{chain}. Given a $v$-chain ${\beta}_1=(r_1,c_1)>\ldots>{\beta}_t=(r_t,c_t)$, we define 
$$s_{\beta_1}\ldots s_{\beta_t}v:=(\{v_1,\ldots,v_d\}\setminus \{c_1,\ldots,c_t\})\cup \{r_1,\ldots ,r_t\}.$$
We say that an element $w$ of $I(d,2d)$ \textbf{dominates} the $v$-chain ${\beta}_1>\ldots>{\beta}_t$ if $w \geq s_{\beta_1}\ldots s_{\beta_t}v $.
Let $\mathfrak{S}$ be a monomial in $\mathfrak{R}(v)$. By a $v$-chain in $\mathfrak{S}$, we mean a sequence ${\beta_1}>\ldots>{\beta_t}$ of elements of $\mathfrak{S}\cap \mathfrak{N}(v)$. We say that $w$ \textbf{dominates} $\mathfrak{S}$ if $w$ dominates every $v$-chain in $\mathfrak{S}$.\par 
We call \textbf{distinguished} the subsets $\mathfrak{S}$ of $\mathfrak{N}(v)$ satisfying the following conditions:

\noindent (A) For $(r,c)\neq (r',c')$ in $\mathfrak{S}$, we have $r\neq r'$, and $c\neq c'$.\\
\noindent (B) If $\mathfrak{S}=\{(r_1,c_1),\ldots,(r_p,c_p)\}$ with $r_1<r_2<\ldots<r_p$, then for $j$, $1\leq j\leq p-1$, we have either $c_j>c_{j+1}$ or $r_j<c_{j+1}$.\par
Example \ref{e.distinguished} below gives an illustration of a distinguished subset. 
\begin{ex}\label{e.distinguished}
{\rm
For $d=7$, $v=(1,2,4,7,8,11,14)$, the subset $\mathfrak{S}$ of $\mathfrak{N}(v)$ given by $$\mathfrak{S}=\{(3,2),(5,1),(9,8),(12,7)\}$$ is distinguished.}
\end{ex}
\begin{rem}
\label{r.distinguished}
By \cite[Proposition 4.3]{Ko-Ra}, there exists a bijection between elements $w$ of $I(d,2d)$ satisfying $w\geq v$ on the one hand, and distinguished subsets of $\mathfrak{N}(v)$ on the other hand. We denote this bijective correspondence by $w\leftrightarrow\mathfrak{S}_w$. 
\end{rem}
Let $\mathfrak{S}$ be a non-empty monomial in $\mathfrak{N}(v)$. If $\beta_1>\cdots>\beta_t$ is a $v$-chain in $\mathfrak{S}$, then we call $\beta_1$ the \textbf{head} of the $v$-chain and $\beta_t$ its \textbf{tail}. We call $t$ to be the \textbf{length} of the $v$-chain. We say that an element $\beta$ of $\mathfrak{S}$ is $t$-\textbf{deep} in $\mathfrak{S}$ (where $t$ is a positive integer) if $\beta$ is the tail of a $v$-chain in $\mathfrak{S}$ of length $t$. The \textbf{depth} of $\beta$ in $\mathfrak{S}$ is defined to be $t$ if $\beta$ is $t$-deep in $\mathfrak{S}$ but not $(t+1)$-deep in $\mathfrak{S}$.\par
Example \ref{e.depth} below gives an illustration of deep and depth of an element in some monomial.
\begin{ex}\label{e.depth}
{\rm
Let $v=(1,2,4,7,8,11,14)$ and $\mathfrak{S}=\{(9,1),(6,2),(5,4),(13,8),(12,11)\}$ be a monomial in $\mathfrak{N}(v)$. Let $\beta=(5,4)$. Then it is easy to see that $\beta$ is $1$-deep, $2$-deep, and $3$-deep in $\mathfrak{S}$. But $\beta$ is not $4$-deep in $\mathfrak{S}$. In fact, $(9,1)>(6,2)>(5,4)$ is a $v$-chain in $\mathfrak{S}$, and this is a $v$-chain in $\mathfrak{S}$ of maximum length having $\beta=(5,4)$ as its tail. Hence the depth of $\beta$ in $\mathfrak{S}$ is $3$ here.}
\end{ex}
We will now recall the map $\pi$ of \cite{Ko-Ra}. Let $\mathfrak{S}$ be a non-empty monomial in the elements of $\mathfrak{N}(v)$. We partition $\mathfrak{S}$ in two stages. First we partition $\mathfrak{S}$ into subsets $\mathfrak{S}_1,\ldots,\mathfrak{S}_k$, where $k$ is the largest length of a $v$-chain in $\mathfrak{S}$: $\beta\in\mathfrak{S}$ belongs to $\mathfrak{S}_j$ if it is $j$-deep but not $(j+1)$-deep.\par
Now we partition each $\mathfrak{S}_j$ into subsets called \textbf{blocks} as follows. We arrange the elements of $\mathfrak{S}_j$ in non-decreasing order of their row numbers (all arrangements are from left to right; and elements occur with their respective multiplicities). Among those with the same row number, the arrangement is by non-decreasing order of column numbers. Two consecutive members $(r,c)$, $(R,C)$ in this arrangement are said to be \textbf{related} if $r>C$. The blocks are the equivalence classes of the smallest equivalence relation containing the above relations.\par
Let $\mathfrak{B}$ be a single block of some $\mathfrak{S}_j$. Let
$$(r_1,c_1),\ldots,(r_p,c_p)$$
be the elements of $\mathfrak{B}$ written in non-decreasing order of both row and column numbers (in such an arrangement, the elements occur with their respective multiplicities). We set $w(\mathfrak{B}):=(r_p,c_1)$ and $\mathfrak{B}'$ to be the monomial
$$\{(r_1,c_2),(r_2,c_3),\ldots,(r_{p-2},c_{p-1}),(r_{p-1},c_p)\}.$$
Set $\mathfrak{S}_j^{(1)}:=\cup_{\mathfrak{B}}\mathfrak{B}'$ (where the index $\mathfrak{B}$ runs over all blocks of $\mathfrak{S}_j$) and $\mathfrak{S}^{(1)}:=\cup_{j=1}^{k}\mathfrak{S}_j^{(1)}$. It follows from \cite[Corollary 4.13]{Ko-Ra} that the set
$$\{w(\mathfrak{B})|\mathfrak{B}\ \text{is a block of}\ \mathfrak{S}\}$$
is a distinguished subset of $\mathfrak{N}(v)$. Let $w$ be the corresponding element of $I(d,2d)$ (under the correspondence given in Remark \ref{r.distinguished}). Set 
$$\pi(\mathfrak{S}):=(w,\mathfrak{S}^{(1)}).$$
This finishes the description of the map $\pi$ of \cite{Ko-Ra}.\par 
Example~\ref{eg.themappi} below gives a detailed illustration of the map $\pi$ of \cite{Ko-Ra}.
\begin{ex}
{\rm
\label{eg.themappi}
Let $d=7$ and $v=(1,2,4,7,8,11,14)$. The dark circles in the grid in Figure~\ref{figure.3} represent a monomial $\mathfrak{S}$ in $\mathfrak{N}^v$, where 
$$\mathfrak{S}=\{(3,2),(5,4),(6,2),(9,1),(9,1),(10,7),(10,7),(10,7),(10,8),(12,1),(13,4)\}.$$
The numbers written near the dark circles denote the multiplicities of these elements in the monomial $\mathfrak{S}$. For this monomial $\mathfrak{S}$, we have
$$\mathfrak{S}_1=\{(9,1),(9,1),(12,1),(13,4)\}$$
$$\mathfrak{S}_2=\{(3,2),(6,2),(10,7),(10,7),(10,7),(10,8)\}$$
$$\mathfrak{S}_3=\{(5,4)\}$$.
Here $\mathfrak{S}_1$ and $\mathfrak{S}_3$ are single blocks. And $\mathfrak{S}_2$ has two blocks given by 
$$\{(3,2),(6,2)\}\ \mbox{and}\ \{(10,7),(10,7),(10,7),(10,8)\}.$$ The dark line segments on the grid show the block decomposition of the monomial $\mathfrak{S}$. The set 
$$\{w(\mathfrak{B})|\mathfrak{B}\ \mbox{is a block of}\ \mathfrak{S}\}=\{(13,1),(6,2),(10,7),(5,4)\}.$$

\vspace{.5 cm}

\begin{figure}[h]
\setlength{\unitlength}{0.6mm}
 \centering
\begin{picture}(60,120)(0,0)
\matrixput(0,0)(15,0){7}(0,15){7}{\circle*{.9}}
\multiputlist(-15,0)(0,15){13,12,10,9,6,5,3}
\multiputlist(0,105)(15,0){1,2,4,7,8,11,14}
\put(0,15){\makebox(0,0){$\bullet$}}
\put(0,45){\makebox(0,0){$\bullet$}}
\put(15,60){\makebox(0,0){$\bullet$}}
\put(15,90){\makebox(0,0){$\bullet$}}
\put(30,0){\makebox(0,0){$\bullet$}}
\put(30,75){\makebox(0,0){$\bullet$}}
\put(45,30){\makebox(0,0){$\bullet$}}
\put(60,30){\makebox(0,0){$\bullet$}}
\linethickness{.7pt}
\dottedline{0.1}(0,45)(0,0)(30,0)
\dottedline{0.1}(45,30)(60,30)
\dottedline{0.1}(15,90)(15,60)
\put(2,17){1}
\put(2,47){2}
\put(17,62){1}
\put(17,92){1}
\put(32,2){1}
\put(32,77){1}
\put(47,32){3}
\put(62,32){1}
\linethickness{.2pt}
\dottedline{0.1}(0,90)(15,90)(15,75)(30,75)(30,45)(60,45)(60,15)(75,15)(75,0)
\put(60,20){$\longleftarrow$}
\put(70,20){ boundary of $\widetilde{\mathfrak{N}(v)}$}
\end{picture}
\caption{Block decomposition of a monomial $\mathfrak{S}$}
\label{figure.3}
\end{figure}
Therefore 
$$w=(5,6,8,10,11,13,14)\ \mbox{and}$$
$$\mathfrak{S}^{(1)}=\{(9,1),(9,1),(12,4),(3,2),(10,7),(10,7),(10,8)\}.$$}
\end{ex}
\bdefn\label{d.SMvvtilde}
Let $v\in I(d,2d)$. A \textbf{standard sequence} in $I(d,2d)$ is a totally ordered sequence $\theta_1\geq\cdots\geq\theta_t$ of elements of $I(d,2d)$. Such a sequence is called $v$-\textbf{compatible} if each $\theta_j$ is comparable to $v$ but no $\theta_j$ equals $v$; it is called \textbf{anti-dominated} by $v$ if $\theta_t\geq v$. Let $\widetilde{SM^{v,v}}$ denote the set of all $v$-compatible standard sequences in $I(d,2d)$ anti-dominated by $v$.
\edefn 
Example \ref{e.SMvvtilde} below gives an illustration of $\widetilde{SM^{v,v}}$.
\begin{ex}\label{e.SMvvtilde}
{\rm
Let $d=4$ and $v=(1,3,5,6)$. Then $(2,5,7,8)\geq (1,4,6,8)\geq (1,4,6,7)$ is an element of $\widetilde{SM^{v,v}}$.}
\end{ex}
Using $\pi$, we now recall the  map $\tilde{\pi}$ of \cite{Ko-Ra} from $mon{\mathfrak{N}(v)}$ to $\widetilde{{SM^{v,v}}}$. Proceed by induction on the degree of an element $\mathfrak{S}$ of $mon{\mathfrak{N}(v)}$. The image of the empty monomial under $\tilde{\pi}$ is taken to be the empty monomial. Let $\mathfrak{S}$ be non-empty, and suppose that $\pi(\mathfrak{S}) =(w,\mathfrak{S}^{(1)})$. By (1) and
(2) of \cite[Proposition~4.1]{Ko-Ra}, the degree of $\mathfrak{S}^{(1)}$ is strictly less than that of $\mathfrak{S}$,  and so by induction $\tilde{\pi}(\mathfrak{S}^{(1)})$ is defined.
Suppose that $\tilde{\pi}(\mathfrak{S}^{(1)})=w'\geq\ldots$. By induction we also know that the degree of $\mathfrak{S}^{(1)}$ is the same as that of $w'\geq\ldots$ and that $w'$ is the least element of $I(d,2d)$ that dominates $\mathfrak{S}^{(1)}$. By (3) of \cite[Proposition~4.1]{Ko-Ra}, we have $w\geq w'$, and we set
$\tilde{\pi}(\mathfrak{S}):=w\geq\tilde{\pi}(\mathfrak{S}^{(1)})$. This finishes the description of the map $\tilde{\pi}$ of \cite{Ko-Ra}.\par
Example \ref{eg.tildepi} below gives an illustration of the map $\tilde{\pi}$ for a monomial in $\mathfrak{N}(v)$.
\begin{ex}\label{eg.tildepi}
{\rm
For the monomial $\mathfrak{S}$ in Example~\ref{eg.themappi} above, we have \par 
\begin{center}
$\tilde{\pi}(\mathfrak{S})=(5,6,8,10,11,13,14)\geq$ $(3,4,8,10,11,12,14)\geq (2,4,7,8,10,11,14)\geq (1,2,7,8,10,11,14)\geq (1,2,4,8,9,11,14)\geq (1,2,4,7,9,11,14)$.
\end{center}}
\end{ex}
\section{\textbf{Statement of the main theorem}}\label{s.main}
\subsection{\textbf{\textit{Extension of the Kodiyalam-Raghavan maps}}}
In this subsection, we will extend the map $\tilde{\pi}$ of Kodiyalam-Raghavan \cite{Ko-Ra} to the entire $\mathfrak{R}(v)$.\par
Fix an element $v$ in $I(d,2d)$. Let $v=(v_1,\ldots,v_d)$. Given $\beta_1=(r_1,c_1)$ and $\beta_2=(r_2,c_2)$ in $\mathfrak{R}(v)\setminus \mathfrak{N}(v)$, we say that $\beta_1>\beta_2$ if $r_1<r_2$ and $c_2<c_1$. A sequence $\beta_1>\cdots>\beta_t$ of elements of $\mathfrak{R}(v)\setminus \mathfrak{N}(v)$ is called an \textbf{anti}-$v$-\textbf{chain}. Given an anti-$v$-chain $\beta_1=(r_1,c_1)>\cdots>\beta_t=(r_t,c_t)$, we define
$$s_{\beta_1}\cdots s_{\beta_t}v:=(\{v_1,\ldots,v_d\}\setminus\{c_1,\ldots,c_t\})\cup\{r_1,\ldots,r_t\}.$$
We say that an element $w$ of $I(d,2d)$ \textbf{anti-dominates} the anti-$v$-chain $\beta_1>\cdots>\beta_t$ if $w\leq s_{\beta_1}\cdots s_{\beta_t}v$. Let $\mathfrak{S}$ be a monomial in $\mathfrak{R}(v)\setminus \mathfrak{N}(v)$. We say that $w$ \textbf{anti-dominates} $\mathfrak{S}$ if $w$ anti-dominates every anti-$v$-chain in $\mathfrak{S}$.\par 
We call \textbf{distinguished} the subsets $\mathfrak{S}$ of $\mathfrak{R}(v)\setminus \mathfrak{N}(v)$ satisfying the following conditions:
\begin{enumerate}
\item  For $(r,c)\neq (r',c')$ in $\mathfrak{S}$, we have $r\neq r'$ and $c\neq c'$.\\
\item If $\mathfrak{S}=\{(r_1,c_1),\ldots,(r_p,c_p)\}$ with $r_1>r_2>\ldots>r_p$, then for $j$, $1\leq j\leq p-1$, we have either $c_j<c_{j+1}$ or $r_j>c_{j+1}$.
\end{enumerate}
\begin{rem}\label{r.anti-distinguished}
It can be proved similarly as in \cite[Proposition 4.3]{Ko-Ra} that there exists a bijection between elements $w$ of $I(d,2d)$ satisfying $w\leq v$ on the one hand and distinguished subsets of $\mathfrak{R}(v)\setminus \mathfrak{N}(v)$ on the other hand. We denote this bijective correspondence by $w\leftrightarrow\mathfrak{S}_w$. 
\end{rem}
Let $\mathfrak{S}$ be a non-empty monomial in $\mathfrak{R}(v)\setminus \mathfrak{N}(v)$. If $\beta_1>\cdots>\beta_t$ is an anti-$v$-chain in $\mathfrak{S}$, then we call $\beta_1$ the \textbf{head} of the anti-$v$-chain and $\beta_t$ its \textbf{tail}. We call $t$ to be the \textbf{length} of the anti-$v$-chain. We say that an element $\beta$ of $\mathfrak{S}$ is $t$-\textbf{deep} in $\mathfrak{S}$ (where $t$ is a positive integer) if $\beta$ is the tail of an anti-$v$-chain in $\mathfrak{S}$ of length $t$. The \textbf{depth} of $\beta$ in $\mathfrak{S}$ is defined to be $t$ if $\beta$ is $t$-deep in $\mathfrak{S}$ but not $(t+1)$-deep in $\mathfrak{S}$.\par
We will now define the map $\pi$ on any monomial in $\mathfrak{R}(v)\setminus \mathfrak{N}(v)$. Let $\mathfrak{S}$ be a non-empty monomial in the elements of $\mathfrak{R}(v)\setminus \mathfrak{N}(v)$. We partition $\mathfrak{S}$ in two stages. First we partition $\mathfrak{S}$ into subsets $\mathfrak{S}_1,\ldots,\mathfrak{S}_k$, where $k$ is the largest length of an anti-$v$-chain in $\mathfrak{S}$: $\beta\in\mathfrak{S}$ belongs to $\mathfrak{S}_j$ if it is $j$-deep but not $(j+1)$-deep.\par
Now we partition each $\mathfrak{S}_j$ into subsets called \textbf{blocks} as follows. We arrange the elements of $\mathfrak{S}_j$ in non-increasing order of their row numbers (where elements occur with their respective multiplicities). Among those with the same row number, the arrangement is by non-increasing order of column numbers. Two consecutive members $(r,c)$, $(R,C)$ in this arrangement are said to be \textbf{related} if $r<C$. The blocks are the equivalence classes of the smallest equivalence relation containing the above relations.\par
Let $\mathfrak{B}$ be a single block of some $\mathfrak{S}_j$. Let
$$(r_1,c_1),\ldots,(r_p,c_p)$$
be the elements of $\mathfrak{B}$ written in non-increasing order of both row and column numbers (in such an arrangement, the elements occur with their respective multiplicities). We set $w(\mathfrak{B}):=(r_p,c_1)$ and $\mathfrak{B}'$ to be the monomial
$$\{(r_1,c_2),(r_2,c_3),\ldots,(r_{p-2},c_{p-1}),(r_{p-1},c_p)\}.$$
Set $\mathfrak{S}_j^{(1)}:=\cup_{\mathfrak{B}}\mathfrak{B}'$ (where the index $\mathfrak{B}$ runs over all blocks of $\mathfrak{S}_j$) and $\mathfrak{S}^{(1)}:=\cup_{j=1}^{k}\mathfrak{S}_j^{(1)}$. It follows (similarly as in \cite[Corollary 4.13]{Ko-Ra}) that the set
$$\{w(\mathfrak{B})|\mathfrak{B}\ \text{is a block of}\ \mathfrak{S}\}$$
is a distinguished subset of $\mathfrak{R}(v)\setminus \mathfrak{N}(v)$. Let $w$ be the corresponding element of $I(d,2d)$ (under the correspondence given in Remark \ref{r.anti-distinguished}). Set 
$$\pi(\mathfrak{S}):=(w,\mathfrak{S}^{(1)}).$$
This finishes the description of the map $\pi$.\par A \textbf{standard sequence} in $I(d,2d)$ is a totally ordered sequence $\theta_1\geq\cdots\geq\theta_t$ of elements of $I(d,2d)$. A standard sequence $\theta_1\geq\cdots\geq\theta_t$ in $I(d,2d)$ is called \textbf{dominated by $v$} if $v\geq\theta_1$. Such a sequence is called $v$-\textbf{compatible} if each $\theta_j$ is comparable to $v$ but no $\theta_j$ equals $v$. Let $\widetilde{SM^v_v}$ denote the set of all \textbf{$v$-compatible standard sequences} in $I(d,2d)$ dominated by $v$.\par
Using $\pi$, we now define the map $\tilde{\pi}$ from the set of all monomials in $\mathfrak{R}(v)\setminus \mathfrak{N}(v)$ to $\widetilde{SM^v_v}$. We proceed by induction on the degree of a monomial $\mathfrak{S}$ in $\mathfrak{R}(v)\setminus \mathfrak{N}(v)$. The image of the empty monomial under $\tilde{\pi}$ is taken to be the empty monomial. Let $\mathfrak{S}$ be non-empty, and suppose that $\pi(\mathfrak{S})=(w,\mathfrak{S}^{(1)})$. It can be shown (similarly as in (1) and (2) of  \cite[Proposition 4.1]{Ko-Ra}) that the degree of $\mathfrak{S}^{(1)}$ is strictly less than that of $\mathfrak{S}$, and so by induction $\tilde{\pi}(\mathfrak{S}^{(1)})$ is defined. Suppose that $\tilde{\pi}(\mathfrak{S}^{(1)})=w'\leq\cdots$. It can be shown (similarly as in (3) of \cite[Proposition 4.1]{Ko-Ra}) that $w\leq w'$. We set $\tilde{\pi}(\mathfrak{S}):=w\leq\tilde{\pi}(\mathfrak{S}^{(1)})$. This finishes the description of the map $\tilde{\pi}$ on the set of all monomials in $\mathfrak{R}(v) \setminus \mathfrak{N}(v)$.\par
Example \ref{eg.tildepifornegative} below gives an illustration of the map $\tilde{\pi}$
for a monomial in $\mathfrak{R}(v) \setminus \mathfrak{N}(v)$.
\begin{ex}\label{eg.tildepifornegative}
{\rm
Let $d=6$ and $v=(3,6,8,10,11,12)$. Let 
$$\mathfrak{S}=\{(9,11),(4,11),(7,10),(5,10),(7,8),(1,8),(4,6)\}$$
be a finite monomial in $\mathfrak{R}(v) \setminus \mathfrak{N}(v)$. 
\begin{figure}[h]
\setlength{\unitlength}{0.8mm}
\centering
\begin{picture}(60,70)(0,0)
\matrixput(0,0)(10,0){5}(0,10){1}{\circle*{.5}}
\matrixput(0,10)(10,0){4}(0,10){2}{\circle*{.5}}
\matrixput(0,30)(10,0){3}(0,10){1}{\circle*{.5}}
\matrixput(0,40)(10,0){2}(0,10){1}{\circle*{.5}}
\linethickness{.5pt}
\dottedline{0.1}(0,50)(0,20)(0,0)(20,0)
\dottedline{0.1}(10,40)(10,20)(30,20)
\put(0,20){\makebox(0,0){$\bullet$}}
\put(0,50){\makebox(0,0){$\bullet$}}
\put(10,30){\makebox(0,0){$\bullet$}}
\put(10,40){\makebox(0,0){$\bullet$}}
\put(20,0){\makebox(0,0){$\bullet$}}
\put(20,40){\makebox(0,0){$\bullet$}}
\put(30,20){\makebox(0,0){$\bullet$}}
\linethickness{0.3pt}
\multiput(-10,0)(1.5,0){34}{\line(1,0){1}}
\multiput(-10,10)(1.5,0){34}{\line(1,0){1}}
\multiput(-10,20)(1.5,0){27}{\line(1,0){1}}
\multiput(-10,30)(1.5,0){27}{\line(1,0){1}}
\multiput(-10,40)(1.5,0){21}{\line(1,0){1}}
\multiput(-10,50)(1.5,0){14}{\line(1,0){1}}
\multiput(-10,0)(0,1.5){34}{\line(0,1){1}}
\multiput(0,0)(0,1.5){34}{\line(0,1){1}}
\multiput(10,0)(0,1.5){34}{\line(0,1){1}}
\multiput(20,0)(0,1.5){27}{\line(0,1){1}}
\multiput(30,0)(0,1.5){20}{\line(0,1){1}}
\multiput(40,0)(0,1.5){7}{\line(0,1){1}}
\multiputlist(-20,0)(0,10){1,2,4,5,7,9}
\put(-13,54){12}
\put(-3,54){11}
\put(8,54){10}
\put(20,43){8}
\put(30,33){6}
\put(40,13){3}
\put(2,15){1}
\put(2,45){1}
\put(12,25){1}
\put(6,35){1}
\put(23,2){1}
\put(15,35){1}
\put(25,15){1}
\end{picture}
\caption{The monomial $\mathfrak{S}$ in $\mathfrak{R}(v) \setminus \mathfrak{N}(v)$ and its block decomposition.}
\label{figure.6}
\end{figure}
Figure~\ref{figure.6} shows the monomial $\mathfrak{S}$ and its block decomposition. The dark circles in the figure represent points in the monomial $\mathfrak{S}$ with their respective multiplicities (the multiplicity of each point in the monomial $\mathfrak{S}$ is $1$ here, which is written near those points in the grid). The dark line segments (together with the point $(7,8)$) denote the blocks of $\mathfrak{S}$.\par
For this monomial $\mathfrak{S}$, we have $\pi(\mathfrak{S})=(w_0,\mathfrak{S}^{(1)})$, where 
$$w_0=(1,3,4,6,7,12)\ \mbox{and}\ \mathfrak{S}^{(1)}=\{(9,11),(4,8),(7,10),(5,6)\}.$$
Then $\pi(\mathfrak{S}^{(1)})=(w_1,\mathfrak{S}^{(2)})$, where $w_1=(3,4,5,8,10,12)$ and $\mathfrak{S}^{(2)}=\{(9,10),(7,8)\}$. And finally, $\pi(\mathfrak{S}^{(2)})=(w_2,\emptyset)$, where $w_2=(3,6,7,9,11,12)$.\par 
Now for the above monomial $\mathfrak{S}$ we have, $$\tilde{\pi}(\mathfrak{S})=(1,3,4,6,7,12) \leq (3,4,5,8,10,12) \leq (3,6,7,9,11,12).$$}
\end{ex}
The relation between the maps $\tilde{\pi}$ and $BRSK$ was given by \cite[Corollary 2.3.2]{papi1}, which is stated here as the following Proposition.
\begin{prop}\label{c.corollary301}
For any monomial $U$ in $\mathfrak{R}(v)$, $\tilde{\pi} = BRSK(U)$.
\end{prop}
\subsection{\textbf{\textit{Ideals of tangent cones to Richardson varieties}}}\label{ss.ideal-of-tgtcone}
Let $\gb$ be the element of $\id$, which was fixed at the beginning of this section. Consider the matrix of size $2d\times d$ whose columns are numbered by the entries of $\gb$, the rows by $\{1,\ldots,2d\}$, the rows corresponding to the entries of $\gb$ form the $d\times d$ identity matrix, and the remaining $d$ rows form a matrix whose entries are $X_{(r,c)}$ such that $(r, c)\in\androotsb$, where $X_{(r,c)}=-X_{(c^*,r^*)}$ if either $r>d$ and $c^*<d$ or $r<d$ and $c^*>d$, and $X_{(r,c)}=X_{(c^*,r^*)}$ otherwise.\par 
For $d=4$, $\beta = (1,2,5,6)$, the $2d \times d$ matrix is given in below:
\begin{center}
\begin{equation}\label{matrix}
\begin{pmatrix}
1 & 0 & 0 & 0\\
0 & 1 & 0 & 0\\
x_{31} & x_{32} & x_{35} & x_{36}\\
x_{41} & x_{42} & x_{45} & x_{35}\\
0 & 0 & 1 & 0\\
0 & 0 & 0 & 1\\
x_{71} & x_{72} & x_{42} & -x_{32}\\
x_{81} & x_{71} & x_{41} & -x_{31}
\end{pmatrix}
\end{equation}
\end{center}
Let  $\miso\subseteq G_{d}(V)\hookrightarrow \mathbb{P}(\bigwedge^d V)$ be the Pl\"ucker embedding (where $G_d(V)$ denotes the
Grassmannian of all $d$-dimensional subspaces of $V$). For $\theta$ in $I(d, 2d)$, let $p_\theta$ denote the corresponding Pl\"ucker coordinate. Consider the affine patch $\affinev$ of $\mathbb{P}(\bigwedge^d V)$ given by $p_\gb\neq 0$. The affine patch $\affinev^\gb:=\miso\cap\affinev$ of the symplectic 
Grassmannian $\miso$ is an affine space whose coordinate ring can be taken
to be the polynomial ring in variables of the form $X_{(r,c)}$ with $(r, c)\in\orootsb$.\par
For $\theta\in I(d,2d)$, consider the submatrix of the above mentioned matrix given by the rows
numbered $\theta\setminus\gb$ and columns numbered $\gb\setminus\theta$. Let $f_{\theta,\gb}$ denote the determinant of this submatrix. Clearly, $f_{\theta,\gb}$ is a homogeneous polynomial in the variables $X_{(r,c)}$, where $(r,c)\in\orootsb$.\par 
Example~\ref{e.fthetabeta} below gives an illustration of $f_{\theta,\beta}$.\par 
\begin{ex}\label{e.fthetabeta}
{\rm
Let $d=4, \ \beta=(1,2,5,6)$, and $\theta =(1,3,4,5)$. So $\theta \in I(4,8)$. Now $f_{\theta,\beta}$ is the determinant of the submatrix whose rows are numbered by $\theta \setminus \beta = \{3,4\}$ and columns are numbered by $\beta \setminus \theta = \{2,6\}$, that is
\begin{equation}
f_{\theta,\beta} = \begin{vmatrix}
x_{32} & x_{36} \\
x_{42} & x_{35}
\end{vmatrix} = x_{32}x_{35}-x_{36}x_{42}.
\end{equation}
Clearly, $f_{\theta,\beta}$ is a homogeneous polynomial in the variables $X_{(r,c)}$, where $(r,c) \in \mathfrak{OR}(\beta)$.}
\end{ex}
The \textbf{$\epsilon$-degree} of an element $x$ of $I(d)$ is defined as the cardinality of $x\setminus [d]$ or equivalently that of $[d]\setminus x$. An ordered pair $\mathfrak{w}=(x,y)$ of elements of $I(d)$ is called an \textbf{admissible pair} if $x\geq y$ and the $\epsilon$-degrees of $x$ and $y$ are equal. We refer to $x$ and $y$ as the \textbf{ top} and the \textbf{bottom} of $\mathfrak{w}$ and write ${\rm top}{\mathfrak{(w)}}$ for $x$ and ${\rm bot}\mathfrak{(w)}$ for $y$. Given any admissible pairs $\mathfrak{w}=(x,y)$ and $\mathfrak{w}'=(x',y')$, we say that $\mathfrak{w}\geq\mathfrak{w}'$ if $y\geq x'$, that is, if $x\geq y\geq x'\geq y'$. 
Let $\mathfrak{w}=(x,y)$ be an admissible pair. Let $\theta$ be the element $(x\cap[d])\cup(y\cap[d]^c)$ of $I(d,2d)$ (as mentioned in \cite[Proposition 3.4]{gr}). For any admissible pair $\mathfrak{w}$, let us denote by $f_{\mathfrak{w},\beta}$ the polynomial $f_{\theta,\beta}$.\par
Example \ref{eg.admipair} below gives an illustration of the admissible pairs.
\begin{ex}\label{eg.admipair}
{\rm
Let $d=4$, so $\epsilon=(1,2,3,4)$. Let $x=(1,4,6,7)$. The $\epsilon$-degree of $x$ is $2$. Let $y=(1,3,5,7)$. Clearly, the $\epsilon$-degree of $y$ is also $2$, and $x,y\in I(d)$ with $x \geq y$. Hence $\mathfrak{w}=(x,y)$ is an admissible pair. Also here ${\rm top}\mathfrak{(w)}=(1,4,6,7)$ and ${\rm bot}\mathfrak{(w)}=(1,3,5,7)$. Again let $x^\prime =(1,3,5,7)$ and $y^\prime = (1,2,4,6)$. Then the $\epsilon$-degrees of both $x^\prime$ and $y^\prime$ are $1$ and $x^\prime,y^\prime\in I(d)$ with $x^\prime \geq y^\prime$. So $\mathfrak{w}'=(x',y')$ is also an admissible pair. As $x \geq y \geq x' \geq y'$, so $\mathfrak{w} \geq \mathfrak{w}'$. Again for the above $\mathfrak{w}$, $\theta =(1,4,5,7)$. So for this $\theta$ and for $\beta=(1,2,5,6)$, (using the $8 \times 4$ matrix of Statement \ref{matrix}) we have,
\begin{center}
    $f_{\mathfrak{w},\beta}=f_{\theta,\beta}=-x_{32}x_{42} - x_{35}x_{72}$.
\end{center}}
\end{ex}
Set $Y_\ga^\gc(\gb):=X_\ga^\gc\cap\affinev^\gb$. From \cite{La-Se-IV} we can deduce a set of generators for the ideal $I_{\alpha,\beta}^{\gamma}$ of functions on $\affinev^\gb$ vanishing on $Y_\ga^\gc(\gb)$. The following equation gives the generators: 
\begin{equation}\label{eq.ideal}
I_{\alpha,\beta}^{\gamma}=\left(f_{\mathfrak{w},\gb}\st \mathfrak{w}=(x,y)\ \text{is an admissible pair}, \ga\not\leq y\ \text{or}\ x\not\leq\gc\right)
\end{equation}
We are interested in the tangent cone to $X_\ga^\gc$ at $e_\gb$ or, what is the
same, the tangent cone to $Y_\ga^\gc(\gb)$ at the origin. Observe that $f_{\mathfrak{w},\gb}$ is a
homogeneous polynomial. Because of this, $Y_\ga^\gc(\gb)$ itself
is a cone and so equal to its tangent cone at the origin. The ideal of the tangent cone to $X_\ga^\gc$ at $e_\gb$ is
therefore the ideal $I_{\ga,\gb}^{\gc}$ in Equation~(\ref{eq.ideal}).
\subsection{\textbf{\textit{Extended $\beta$-chains}}}\label{ss.extd-beta-chains}
Let $\gb$ be the element of $\id$, which was fixed at the end of \S\ref{ss.In-state-problem}. For elements $\lambda=(R,C),\mu=(r,c)$ of $\rootsb$, we write $\lambda>\mu$, if $R>r$ and $C<c$ (note that these are strict inequalities). A sequence $\lambda_1>\cdots>\lambda_k$ of elements of $\rootsb$ is called an \textbf{extended} $\beta$-\textbf{chain}. Note that an extended $\gb$-chain can also be empty. Letting $C$ to be an extended $\gb$-chain, we define $C^+:=C\cap\posb$ and $C^-:=C\cap(\rootsb\setminus\posb)$. We call $C^+$ (resp. $C^-$) the \textbf{positive} (resp. \textbf{negative}) \textbf{part} of the extended $\gb$-chain $C$. We call an extended $\gb$-chain $C$ \textbf{positive} (resp. \textbf{negative}) if $C=C^+$ (resp. $C=C^-$). The extended $\gb$-chain $C$ is called \textbf{non-vanishing} if at least one of its positive or negative part is non-empty. Clearly then, every non-empty extended $\gb$-chain is non-vanishing.\par 
An extended $\gb$-chain that lies completely in $\orootsb$ is called an \textbf{extended upper} $\gb$-\textbf{chain}. We similarly define \textbf{extended upper positive} and \textbf{extended upper negative} $\gb$-chains.\par 
Example~\ref{e.extendedchain} below illustrates an extended $\gb$-chain and an extended upper $\gb$-chain.\par 
\begin{ex}\label{e.extendedchain}
{\rm
Let $d=7$ and $\beta=(1,3,4,7,9,10,13)$. Let $\lambda_1=(14,1)$, $\lambda_2=(12,3)$, $\lambda_3=(6,7)$, and $\lambda_4=(5,13)$. Then clearly $\lambda_1 > \lambda_2 > \lambda_3>\lambda_4$. So this is an extended $\beta$-chain in $\mathfrak{R}(\beta)$. If we denote the above $\beta$-chain by $C$, then $C^+=\lambda_1 > \lambda_2$ and $C^-=\lambda_3 > \lambda_4$. Again $\lambda_1 > \lambda_2 > \lambda_3$ is an extended upper $\beta$-chain.}
\end{ex}
\bdefn
Let $\gb$ be as fixed earlier. Let $\bar{\gb}$ denote the set $[2d]\setminus\gb$. We call $\bar{\gb}$ the \textbf{complement} of $\gb$.
\edefn
\bdefn\label{d.AminusB}
Let $A\subset\bar{\gb}$ and $B\subset\gb$. We define $A-B$ as the set $A\cup(\gb\setminus B)$.
\edefn
\bdefn\label{d.good-chains}
Let $C$ be an extended upper $\gb$-chain. Let $(P^C,Q^C)=BRSK(C\cup C^{\hash})$. Let $(P_1^C,Q_1^C)$ denote the topmost row of $(P^C,Q^C)$ and $(P_r^C,Q_r^C)$ denote the bottom-most row of $(P^C,Q^C)$. Let ${\rm top}(C^+)$ denote the element $P_r^C-Q_r^C$ of $I(d,2d)$ and ${\rm bot}(C^-)$ denote the element $P_1^C-Q_1^C$ of $I(d,2d)$, where $P_r^C-Q_r^C$ and $P_1^C-Q_1^C$ are having the meaning as given in Definition \ref{d.AminusB}.
\edefn
\begin{thm}\label{t.belongtoId}
Given any extended upper $\gb$-chain $C$, the elements ${\rm top}(C^+)$ and ${\rm bot}(C^-)$ of $I(d,2d)$ in fact belong to $I(d)$. 
\end{thm}
\begin{proof}Let $(C\cup C^{\hash})^+$ and $(C\cup C^{\hash})^-$ denote the positive and negative parts respectively of the multiset $C\cup C^{\hash}$. We know that $BRSK(C\cup C^{\hash})$ is equal to the notched bitableau obtained by placing the notched bitableau $BRSK((C\cup C^{\hash})^-)$ on top of the notched bitableau $BRSK((C\cup C^{\hash})^+)$.

Recall the map $\tilde{\pi}$ from \cite[\S 4]{Ko-Ra}. We know from \cite[Corollary 2.3.2]{papi1} that $BRSK((C\cup C^{\hash})^+)=\tilde{\pi}((C\cup C^{\hash})^+)$. Also $(C\cup C^{\hash})^+=((C\cup C^{\hash})^+)^{\hash}$. It hence follows from \cite[Proposition 5.6]{gr} that all the elements of $I(d,2d)$ corresponding to all the rows of $BRSK((C\cup C^{\hash})^+)$ in fact belong to $I(d)$. In particular, the element $P^C_r-Q^C_r={\rm top}(C^+)$ also belongs to $I(d)$. The proof of the fact that ${\rm bot}(C^-)$ belongs to $I(d)$ is similar [We omit the proof here because it involves proving that the maps $BRSK$ and $\tilde{\pi}$ are equal on negative multisets. And this proof is similar to that in \cite{papi1}]. 
\end{proof}
The example below illustrates Theorem \ref{t.belongtoId}.
\begin{ex}
{\rm
\newcommand{\A}{13}
\newcommand{\B}{10}
\newcommand{\M}{11}
\newcommand{\N}{12}
\newcommand{\K}{14}
Let $d=7$ and $\beta=(1,3,4,7,9,10,13)$. Clearly, $\beta \in I(d)$. Now $\bar{\beta}=(2,5,6,8,11,12,14)$. Consider the upper extended $\beta$-chain 
$$C=\{(12,1),(11,3),(8,4),(6,7),(5,9),(2,10)\}.$$
According to \cite{Kr-bkrs},
\begin{equation}
BRSK(C \cup C^\hash)= \left( \hspace{.2cm}\young(256,568,8\M\N,\M\N\K)\hspace{.2cm} ,  \hspace{.2cm} \young(9\B\A,79\B,347,134) \right) 
\end{equation}
Now ${\rm top}(C^+)=P_r^C - Q_r^C =\{11,12,14\}\cup (\beta \smallsetminus \{1,3,4\})=(7,9,10,11,12,13,14) $, and ${\rm bot}(C^-)=P_1^C - Q_1^C =\{2,5,6\}\cup (\beta \smallsetminus \{9,10,13\})=(1,2,3,4,5,6,7) $. Clearly, both of ${\rm top}(C^+)$ and ${\rm bot}(C^-)$ belong to $I(d)$.}
\end{ex}

\subsection{\textbf{\textit{Gr$\ddot{o}$bner basis for ideals of tangent cones}}}\label{ss.termorder}
We now specify the \textbf{term order $\torder$} on monomials in the coordinate functions $\{X_{(r,c)}|(r,c)\in\orootsb\}$ with respect to which the
initial ideal of the ideal $I_{\ga,\gb}^\gc$ of the tangent cone is to be taken.\par
\bdefn\label{d.term-order}
Let $>$ be the total order on~$\orootsb$ satisfying the following condition:
\begin{itemize}
\item 
$X_{(r,c)}>X_{(r',c')}$ if either (a) $r>r'$ or (b) $r=r'$ and $c<c'$.
\end{itemize}
Let~$\torder$ be the term order on monomials in~$\orootsb$ given by the homogeneous lexicographic order with respect to~$>$. 
\edefn
Example \ref{eg.termorder} below gives an illustration of the term order $\torder$.
\begin{ex}\label{eg.termorder}
{\rm
Let $d=7$ and $\beta=(1,3,4,7,9,10,13)$. Now all of $(14,1),(12,3),(11,3),(11,1)$ are elements in~$\orootsb$, and according to the above term order we have, $X_{(14,1)}>X_{(12,3)}>X_{(11,1)}>X_{(11,3)}$. Let $X_{\mathfrak{S}_1}=X_{(14,1)}^3X_{(11,1)}X_{(11,3)}^2$, $X_{\mathfrak{S}_2}=X_{(14,1)}X_{(12,3)}X_{(11,3)}$, and  $X_{\mathfrak{S}_3}=X_{(14,1)}X_{(11,1)}^2$. Clearly $\mathfrak{S}_1,\ \mathfrak{S}_2$, and $\mathfrak{S}_3$ all are monomials in~$\orootsb$. Now the degree of the polynomial $X_{\mathfrak{S}_1}$ is greater than that of $X_{\mathfrak{S}_2}$ and $X_{\mathfrak{S}_3}$. So $X_{\mathfrak{S}_1} \torder X_{\mathfrak{S}_2}$ and $X_{\mathfrak{S}_1}\torder X_{\mathfrak{S}_3}$. Though the degree of $X_{\mathfrak{S}_2}$ is equal to the degree of $X_{\mathfrak{S}_3}$, but $X_{(12,3)} >X_{(11,1)}$ and in $X_{\mathfrak{S}_2}$, the degree of $X_{(12,3)}$ is one and in $X_{\mathfrak{S}_3}$, the degree of $X_{(12,3)}$ is zero. So according to the definition of homogeneous lexicographic order, we have $X_{\mathfrak{S}_2} \torder X_{\mathfrak{S}_3}$. Hence $X_{\mathfrak{S}_1} \torder X_{\mathfrak{S}_2} \torder X_{\mathfrak{S}_3}$.} 
\end{ex}
Now recall that the ideal of the tangent cone to $X_\ga^\gc$ at $e_\gb$ is
the ideal $I_{\ga,\gb}^\gc$ given by Equation~\ref{eq.ideal}. Let $\torder$ be as in \S \ref{ss.termorder}. For any element $f\in I_{\ga,\gb}^\gc$, let $\init f$ denote the initial term of $f$ with respect to the term order $\torder$. We define $\init I_{\ga,\gb}^\gc$ to be the ideal $\langle\init f\st f\in I_{\ga,\gb}^\gc\rangle$ inside the polynomial ring $P:=K[X_{(r,c)}\st (r,c)\in\orootsb]$. 
\bdefn\label{d.good-pair}
An admissible pair $\mathfrak{w}=(t,u)$ (where $t\geq u$) is called a \textbf{good admissible pair} if it satisfies both of the following 2 properties:
\begin{enumerate}
\item $\ga\nleq u$ or $t\nleq\gc$.
\item Either $\init f_{\mathfrak{w},\gb}$ forms a positive upper extended $\gb$-chain $C^+$ such that $C_{(1)}^{+}-C_{(2)}^{+}\nleq\gc$ or $\init f_{\mathfrak{w},\gb}$ forms a negative upper extended $\gb$-chain $C^-$ such that $C_{(1)}^{-}-C_{(2)}^{-}\ngeq\ga$.
\end{enumerate}
Let $\mathcal{G}_{\ga,\gb}^{\gc}$ denote the set $\{f_{\mathfrak{w},\gb}|\mathfrak{w}\ \text{is good}\}$.
\edefn
Example~\ref{ex.1} below illustrates a good admissible pair.
\begin{ex}\label{ex.1}
{\rm
Let $d=4$, $\alpha=(1,2,3,5),\ \beta=(1,2,5,6)$, and $\gamma=(2,3,5,8)$. Let $\mathfrak{w}=(t,u)$ be an admissible pair, where $t=(3,4,7,8)$ and $u=(1,2,5,6)$. 
Clearly $t \nleq \gamma$. Now in \S~\ref{ss.ideal-of-tgtcone}, we have already defined that $\theta = (t \cap [d]) \cup (u \cap [d]^c)$, and $f_{\mathfrak{w},\beta} = f_{\theta, \beta}$. Hence in this example $\theta = (3,4,5,6)$ and from the matrix which is given in \S~\ref{ss.ideal-of-tgtcone}, we have  
\begin{equation}
f_{\mathfrak{w},\beta} = \begin{vmatrix}
x_{31} & x_{32} \\
x_{41} & x_{42}
\end{vmatrix}
\end{equation}
Observe that $\init f_{\mathfrak{w},\beta}= -x_{41}x_{32}$. Clearly, $\init f_{\mathfrak{w},\gb}$ forms a positive upper extended $\gb$-chain $C^+$ such that $C_{(1)}^{+}-C_{(2)}^{+}= \{3,4\}\cup (\beta \smallsetminus \{1,2\})=(3,4,5,6)\nleq\gc$. Hence $\mathfrak{w}=(t,u)$ is a good admissible pair.}
\end{ex}
\bdefn\label{d.initial}
If $S$ is any nonempty subset of the polynomial ring $P:=K[X_{(r,c)}\st (r,c)\in\orootsb]$ such that $S \neq \{0\}$. We define $\init S$ to be the ideal $\langle \init(s)|s\in S\rangle$.
\edefn
\noindent The main result of this paper is the following:
\begin{thm}\label{t.main}
The set $\mathcal{G}_{\ga,\gb}^{\gc}$ is a Gr\"obner basis for the ideal $I_{\ga,\gb}^\gc$.
\end{thm}
\subsection{\textbf{\textit{Strategy of the proof}}}\label{ss.mainthm+strategy}
To explain the strategy of the proof of Theorem~\ref{t.main}, we need the following definition.
\bdefn\label{d.standard-monomials}
We call $f=f_{\mathfrak{w}_1,\gb}\cdots f_{\mathfrak{w}_r,\gb}\in P=K[X_{(r,c)}|(r,c)\in\orootsb]$ a \textbf{standard monomial} if 
\begin{equation}\label{e.standard1}
\mathfrak{w}_1\leq\cdots\leq\mathfrak{w}_r, 
\end{equation}
and for each $i\in\{1,\ldots,r\}$, we have 
\begin{equation}\label{e.standard2}
\text{Either}\ \gb\geq {\rm top}(\mathfrak{w}_i)\ \text{or}\ {\rm top}(\mathfrak{w}_i)\geq\gb,
\end{equation}
\begin{equation}
\label{e.standard3}
\text{and}\ \text{either}\ \gb\geq \rm {bot}(\mathfrak{w}_i)\ \text{or}\ {\rm bot}(\mathfrak{w}_i)\geq\gb,
\end{equation}
\begin{equation}
\label{e.standard4}
\text{and}\ \mathfrak{w}_i\neq(\gb,\gb). 
\end{equation}
If in addition, for $\ga,\gc\in I(d)$, we have
\begin{equation}\label{e.standard5}
\ga\leq {\rm bot}(\mathfrak{w}_1)\ \text{and}\ {\rm top}(\mathfrak{w}_r)\leq\gc, 
\end{equation}
then we say that $f$ is \textbf{standard on} $Y_{\ga}^{\gc}(\gb)$.
\edefn
Example \ref{eg.standardmono} below gives an illustration of a standard monomial on $Y_{\ga}^{\gc}(\gb)$.
\begin{ex}\label{eg.standardmono}
{\rm
Let $d=4$, $\ga =(1,2,3,5)$, $\gb=(1,2,5,6)$, and $\gc=(3,4,7,8)$. For this $\gb$, the $8\times 4$ matrix is given by \ref{matrix}. Let $\mathfrak{w}_1=((1,2,4,6),(1,2,3,5))$ and $\mathfrak{w}_2=((2,4,6,8),(2,3,5,8))$. Clearly, $\mathfrak{w}_1$ and $\mathfrak{w}_2$ both are admissible pairs. Let $\theta_1$ and $\theta_2$ be the images of $\mathfrak{w}_1$ and $\mathfrak{w}_2$ respectively, under the correspondence given by $\mathfrak{w}=(x,y)\mapsto\theta=(x \cap [d]) \cup (y \cap [d]^c)$ as mentioned in \cite[Proposition 3.4]{gr}. So we have, $\theta_1=(1,2,4,5)$ and $\theta_2=(2,4,5,8)$. Now using the $8 \times 4$ matrix of (\ref{matrix}) we have, 
\begin{center}
    $f_{\mathfrak{w}_1,\beta}=f_{\theta_1,\beta}=x_{35} \in P$\\
    and $f_{\mathfrak{w}_2,\beta}=f_{\theta_2,\beta}=-x_{41}x_{31}-x_{35}x_{81} \in P$.
\end{center}
Hence $f =f_{\mathfrak{w}_1,\beta}f_{\mathfrak{w}_2,\beta} \in P$. Clearly, $\mathfrak{w}_1 \leq \mathfrak{w}_2$. Again $\beta \geq {\rm top}\mathfrak{(w_1)}$ and $\beta \geq {\rm bot}\mathfrak{(w_1)}$. Also $\beta \leq {\rm top}\mathfrak{(w_2)}$ and $\beta \leq {\rm bot}\mathfrak{(w_2)}$, and $\mathfrak{w}_i \neq (\beta,\beta) \ \mbox{ for all } i \in \{1,2\}$. Again $\ga$, $\gc \in I(d)$ are such that $\ga \leq {\rm bot}\mathfrak{(w_1)}$ and $\gc \geq {\rm top}\mathfrak{(w_2)}$. Hence $f$ is standard on $Y_\ga^\gc(\gb)$.}
\end{ex}
\bdefn\label{d.degree-standrad-monomial}
Let $f=f_{\mathfrak{w}_1,\gb}\cdots f_{\mathfrak{w}_r,\gb}$ be a standard monomial on $Y_{\ga}^{\gc}(\gb)$. We define the \textit{degree} of $f$ to be the sum of the $\gb$-degrees of $\mathfrak{w}_1,\ldots,\mathfrak{w}_r$, where given any admissible pair $\mathfrak{w}=(x,y)$, the $\gb$-degree of $\mathfrak{w}$ is defined to be $\frac{1}{2}(|x\setminus\gb|+|y\setminus\gb|)$.
\edefn
We now briefly sketch the proof of Theorem~\ref{t.main} (the details are found in \S \ref{s.theproof}). Clearly, $\mathcal{G}_{\ga,\gb}^{\gc}$ is contained in the ideal $I_{\ga,\gb}^{\gc}$. So $\init \mathcal{G}_{\ga,\gb}^{\gc}\subseteq\init I_{\ga,\gb}^\gc$. Hence to prove Theorem \ref{t.main}, we only need to show that in any degree, the number of monomials of $\init \mathcal{G}_{\ga,\gb}^{\gc}$ is at least as great as the number of monomials of $\init I_{\ga,\gb}^\gc$ (the other inequality being trivial). Equivalently, we need to prove that in any degree, the number of monomials of $P\setminus\init \mathcal{G}_{\ga,\gb}^{\gc}$ is no greater than the number of monomials of $P\setminus\init I_{\ga,\gb}^\gc$. Both the monomials of $P\setminus\init I_{\ga,\gb}^\gc$ and the standard monomials on $Y_{\ga}^{\gc}(\gb)$ (the definition of a standard monomial on $Y_{\ga}^{\gc}(\gb)$ is given in Definition~\ref{d.standard-monomials}) form a basis for $P/I_{\ga,\gb}^\gc$, and thus agree in cardinality in any degree. Therefore it suffices to prove that, in any degree, the number of monomials of $P\setminus\init\mathcal{G}_{\ga,\gb}^{\gc}$ is less than or equal to the number of standard monomials on $Y_{\ga}^{\gc}(\gb)$. In this paper, we consider two sets, namely, the set of all ``non-vanishing special multisets on $\BBZ$ (bounded by $T_{\ga}$, $W_{\gc}$)'', and the set of all ``non-vanishing semistandard notched bitableaux on $(\BBZ)^*$ (bounded by $T_{\ga}$, $W_{\gc}$)''. The meaning attached to these two sets is given in \S \ref{s.two-sets} below. In \S \ref{s.theproof} below, we will first show that there exists a degree doubling injection from the set of all monomials of $P\setminus\init \mathcal{G}_{\ga,\gb}^{\gc}$ to the former set. Then we will show that, there exists a degree-halving injection from the later set (namely, the set of all ``non-vanishing semistandard notched bitableaux on $(\BBZ)^*$ (bounded by $T_{\ga},W_{\gc}$)'') to the set of all standard monomials on $Y_{\ga}^{\gc}(\gb)$.  And then we will prove that the map $BRSK$ of \cite{Kr-bkrs} is a degree preserving bijection from the former set to the later. This will complete the proof.\par
Example \ref{e.grobner-basis} below gives an illustration of a Gr$\ddot{o}$bner basis.
\begin{ex}\label{e.grobner-basis}
{\rm
Let $d=4$, $\ga=(1,2,3,5)$, $\gb=(1,2,5,6)$, and $\gc=(2,3,5,8)$.  Then from Example \ref{ex.1}, we know that $\mathfrak{w}=(t,u)$ is an good admissible pair, where $t=(3,4,7,8)$ and $u=(1,2,5,6)$. For the above $\alpha, \beta, \gamma$ if we consider all the admissible pairs which satisfying both the conditions of good admissible pair, then we will get the set $\mathcal{G}_{\ga,\gb}^{\gc}$ which is given by $\{f_{\mathfrak{w},\gb}|\mathfrak{w}\in G\}$, where $G$ is the following set (of all good admissible pairs):
$$G=\{((1,2,3,4),(1,2,3,4)),((1,4,6,7),(1,2,5,6)), ((2,4,6,8),(1,2,5,6)), ((3,4,7,8),(1,2,5,6)),$$
$$((1,5,6,7),(1,5,6,7)), ((2,5,6,8),(1,5,6,7)), ((3,5,7,8),(1,5,6,7)), ((4,6,7,8),(1,5,6,7)),$$
$$((2,5,6,8),(2,5,6,8)), ((3,5,7,8),(2,5,6,8)), ((4,6,7,8),(2,5,6,8)), ((5,6,7,8),(5,6,7,8))\}.$$
As in Example~\ref{ex.1}, we can easily find the initial term of the above good admissible pairs. In this case, we have
$$\init \mathcal{G}_{\ga,\gb}^{\gc}=\langle\{x_{71}x_{32}, x_{71}, x_{72}, x_{41}x_{32}, x_{41}, x_{42}, x_{45}x_{36}, x_{81}x_{72}, x_{81}x_{42}, x_{81}x_{32}, x_{81}, x_{71}x_{42}\}\rangle.$$}
\end{ex}
\section{\textbf{The two sets}}\label{s.two-sets}
As mentioned towards the end of the previous section, the two sets under consideration are ``non-vanishing special multisets on $\BBZ$ (bounded by $T_{\ga}$, $W_{\gc}$)'' and ``non-vanishing semistandard notched bitableaux on $(\BBZ)^*$ (bounded by $T_{\ga}$, $W_{\gc}$)''. We will now explain the meaning of these two sets.\par
Let $\ga,\gb,\gc$ be as before (\S \ref{ss.In-state-problem}). Let $I_{\gb}$ be the set of all pairs $(R,S)$ such that $R\subset\bar{\gb}$, $S\subset\gb$, and $|R|=|S|$. Let $I_{\gb}^{*}$ be the subset of $I_{\gb}$ consisting of all pairs $(R,S)$ such that $R=S^*$. Clearly then, the map $(R,S)\mapsto R-S$ is a bijection from $I_{\gb}^{*}$ to $I(d)$ (Indeed, the inverse map is given by $\theta\mapsto(\theta\setminus\gb,\gb\setminus\theta)$). Let $(R_{\ga},S_{\ga})$ and $(R_{\gc},S_{\gc})$ be the preimages of $\ga$ and $\gc$ respectively under the bijection from $I_{\gb}^{*}$ to $I(d)$. Define $T_{\ga}$ and $W_{\gc}$ to be any subsets of $\BBZ$ such that $(T_{\ga})_{(1)}=R_{\ga},(T_{\ga})_{(2)}=S_{\ga}, (W_{\gc})_{(1)}=R_{\gc}$, and $(W_{\gc})_{(2)}=S_{\gc}$. Note that there always exist subsets $T_{\ga}$ and $W_{\gc}$ of $\BBZ$ such that $T_{\ga}$ is negative and $W_{\gc}$ is positive [This is because $\gb\leq\gc$. Apply the first half of the proof of \cite[Proposition 4.3]{Ko-Ra} to $\gc$ (which is $\geq\gb$) to get a distinguished monomial corresponding to $\gc$. This distinguished monomial will serve as a positive subset $W_{\gc}$ of $\BBZ$. Similarly, we can get a negative subset $T_{\ga}$ of $\BBZ$ corresponding to $\ga$ (which is $\leq\gb$)]. Hence we can choose $T_{\ga}$ and $W_{\gc}$ in such a way that the former is negative and the later is positive.\par
Example \ref{eg.illustration} below gives an illustration of the above paragraph.
\begin{ex}\label{eg.illustration}
{\rm
Let $d=7$ and $\beta = (1,3,4,7,9,10,13)$. So $\bar{\beta} = (2,5,6,8,11,12,14)$. Let $\alpha = (1,2,3,5,7,9,11)$ and $\gamma = (4,5,6,7,12,13,14)$. So $\alpha \leq \beta \leq \gamma$. Now $I_{\beta}$ is the set of all pairs $(R,S)$ such that $R \subset (2,5,6,8,11.12.14)$, $S \subset (1,3,4,7,9,10,13)$, and $|R| =|S|$. Let $R=(2,6,8,11)$ and $S=(4,7,9,13)$. Clearly, $|R|=|S|$ and $R= S^{\star}$. So according to the definition of $I_{\beta}^{\star}$, $(R,S) \in I_{\beta}^{\star}$. Now $R-S = (1,2,3,6,8,10,11)$ is in $I(d)$. Again both of 
$$(R_{\alpha}, S_{\alpha})=((2,5,11),(4,10,13))\ \mbox{and}\ (R_{\gamma},S_{\gamma})=((5,6,12,14),(1,3,9,10))$$
are in $I_{\beta}^{\star}$. Let $T_{\alpha}=\{(2,4),(5,10),(11,13)\}$ and $W_{\gamma}=\{(5,1),(6,3),(12,9),(14,10)\}$. Then clearly $T_{\alpha}$ is a negative and $W_{\gamma}$ is a positive subset of $\bar{\beta} \times \beta$.}
\end{ex}
\subsection{\textbf{\textit{The first set}}}\label{ss.thefirstset}
A non-vanishing multiset on $\BBZ$ (bounded by $T_{\ga}$, $W_{\gc}$) has the same meaning as in \S\ref{ss.basicnot}. Such a multiset $\mathfrak{S}$ is called a \textbf{non-vanishing special multiset on} $\BBZ$ \textbf{(bounded by $T_{\ga}, W_{\gc}$)} if moreover, the following two properties are satisfied:
\begin{enumerate}
\item $\mathfrak{S}=\mathfrak{S}^{\hash}$.
 \item the multiplicity of any diagonal element in $\mathfrak{S}$ is even.
 \end{enumerate}
Example \ref{e.thefirstset} below gives an illustration of the first set. 
\begin{ex}\label{e.thefirstset}
{\rm
Let $d=7$ and $\beta=(1,3,4,7,9,10,13)$. Let $\alpha=(1,2,3,5,7,9,11)$ and $\gamma=(4,5,6,7,12,13,14)$. Let $\mathfrak{S}=\{(2,3),(12,13),(5,10),(5,10)\}$. Clearly, $\mathfrak{S}=\mathfrak{S}^{\hash}$ and the multiplicity of any diagonal element in $\mathfrak{S}$ is even. The only $\beta$-chains in $\mathfrak{S}$ are $C_1=\{(2,3)\}$, $C_2=\{(12,13)\}$, and $C_3=\{(5,10)\}$. \\
Let us take $T_{\alpha}=\{(2,4),(5,10),(11,13)\}$ and $W_{\gamma}=\{(5,1),(6,3),(12,9),(14,10)\}$. Clearly, $(T_{\alpha})_1-(T_{\alpha})_2=\alpha$ and $(W_{\gamma})_1-(W_{\gamma})_2=\gamma$. We have to check that $T_{\alpha}\leq C_i\leq W_{\gamma}$ for all $i\in\{1,2,3\}$.\\
Now, 
$$\{2\}-\{3\}=\{2\}\cup(\beta\setminus\{3\})=(1,2,4,7,9,10,13).$$
So $T_{\alpha}\leq C_1\leq W_{\gamma}$. Similarly one can check that $T_{\alpha}\leq C_2\leq W_{\gamma}$ and $T_{\alpha}\leq C_3\leq W_{\gamma}$. Hence $\mathfrak{S}$ is a special multiset.}
\end{ex}
\subsection{\textbf{\textit{The second set}}}\label{ss.thesecondset}
A non-vanishing semistandard notched bitableau on $\BBZ$ bounded by $T_{\ga}$, $W_{\gc}$ has the same meaning as in \S\ref{ss.BRSK}. Such a notched bitableau $(P,Q)$ is said to be a \textbf{non-vanishing semistandard notched bitableau on} $(\BBZ)^*$ \textbf{(bounded by $T_{\ga}, W_{\gc}$)} if moreover, the following 5 conditions are satisfied:
\begin{enumerate}
\item  $P_i=Q_i^*$ for every row number $i$ of $(P,Q)$.
\item $(P,Q)$ doesn't contain any empty rows.
\item The total number of rows in $P$ (or $Q$) is either even, or it is odd but
$$P_1-Q_1\leq\cdots\leq P_n-Q_n\leq\beta\leq P_{n+1}-Q_{n+1}\leq\cdots\leq P_{n+p}-Q_{n+p},$$
where $n+p$ is the total number of rows in $P$ (or $Q$), and $(P_i,Q_i)$ (for $1\leq i\leq n$) is the negative part of $(P,Q)$, and $(P_{n+i},Q_{n+i})$ (for $1\leq i\leq p$) is the positive part of $(P,Q)$.\\
Let us denote  by $\delta_1\leq\cdots\leq\delta_{n+p+1}$ the sequence $P_1-Q_1\leq\cdots\leq P_n-Q_n\leq\beta\leq P_{n+1}-Q_{n+1}\leq\cdots\leq P_{n+p}-Q_{n+p}$, where $n+p$ is odd.
\item Either the total number of rows of $P$ (or $Q$) is even, and the $\epsilon$-degrees (where $\epsilon=(1,2,\ldots,d)\in I(d)$) of $P_j-Q_j$ and $P_{j+1}-Q_{j+1}$ are equal for each $j$ odd, or the total number of rows in $P$ (or $Q$) is odd (say, $n+p$), and the $\epsilon$-degrees of $\delta_j$ and $\delta_{j+1}$ are equal for each $j$ odd, where the $\delta_j$'s are as mentioned in item (3) above.
\item The total number of boxes in $P$ (or $Q$) is even.
\end{enumerate}
Example \ref{e.thesecondset} below gives an illustration of the second set.
\begin{ex}\label{e.thesecondset}
{\rm
\newcommand{\A}{10}
\newcommand{\B}{11}
\newcommand{\C}{12}
\newcommand{\D}{13}
\newcommand{\E}{14}
Let $d=7$ and $\beta=(1,3,4,7,9,10,13)$. Let
$$(P,Q)=\left(\hspace{.2cm}\young(2\B,5\C,6\E)\hspace{.2cm},\hspace{.2cm}\young(4\D,3\A,19)\hspace{.2cm}\right).$$
Clearly, $(P,Q)$ is a notched bitableau on $\BBZ$. Let $\alpha=(1,2,3,5,7,9,11)$ and $\gamma=(4,5,6,7,12,13,14)$. Let us take $T_{\alpha}=\{(2,4),(5,10),(11,13)\}$ and $W_{\gamma}=\{(5,1),(6,3),(12,9),(14,10)\}$. Observe that
$$P_1-Q_1=(1,2,3,7,9,10,11),$$
$$P_2-Q_2=(1,4,5,7,9,12,13),$$
$$\mbox{and}\ P_3-Q_3=(3,4,6,7,10,13,14).$$
Since $P_1\ltS Q_1$, $P_2\gtS Q_2$, and $P_3\gtS Q_3$, we have $(P,Q)$ is non-vanishing. Also, $P_1-Q_1\leq P_2-Q_2\leq P_3-Q_3$. Hence $(P,Q)$ is semistandard. Again,
$$(T_{\alpha})_{(1)}-(T_{\alpha})_{(2)}=\alpha\leq P_1-Q_1\ \mbox{and}$$
$$P_3-Q_3\leq \gamma=(W_{\gamma})_{(1)}-(W_{\gamma})_{(2)}.$$
So $(P,Q)$ is bounded by $T_{\alpha}$ and $W_{\gamma}$. Observe now that
\begin{enumerate}
\item $P_i=Q_i^*$ for all $i\in\{1,2,3\}$.
\item $(P,Q)$ does not contain any empty rows.
\item The total number of rows in $P$ (or $Q$) is $3$, which is odd, but
$$P_1-Q_1\leq\beta\leq P_2-Q_2\leq P_3-Q_3.$$
\item The $\epsilon$-degrees of $P_1-Q_1$ and $\beta$ are the same (both are $3$). Also, the $\epsilon$-degrees of $P_2-Q_2$ and $P_3-Q_3$ are the same (both are $3$).
\item The total number of boxes in $P$ (or $Q$) is $6$, which is even.
\end{enumerate}
Hence $(P,Q)$ is on $(\BBZ)^*$.}
\end{ex}
\section{\textbf{The proof}}\label{s.theproof}
The main result (Theorem \ref{t.main}) will be obtained as a consequence of Theorem \ref{t.the-first-set}, \ref{t.the-second-set}, and \ref{t.brsk-bijection}. For this section, we fix $\alpha, \beta$, and $\gamma$ in $I(d)$ such that $\alpha \leq \beta \leq \gamma$. Also we need some definitions for this section, which we will state first.
\bdefn\label{d.SMvv}
An ordered sequence $(\mathfrak{w}_1,\ldots,\mathfrak{w}_t)$ of admissible pairs is called a \textbf{standard sequence of admissible pairs} if $\mathfrak{w}_i\geq\mathfrak{w}_{i+1}$ for $1\leq i<t$. We often write $\mathfrak{w}_1\geq\cdots\geq\mathfrak{w}_t$ to denote the standard sequence $(\mathfrak{w}_1,\ldots,\mathfrak{w}_t)$ of admissible pairs.
Given any $v\in I(d)$, we say that a standard sequence $\mathfrak{w}_1\geq\cdots\geq\mathfrak{w}_t$ of admissible pairs is $v$-\textbf{compatible} if for each $\mathfrak{w}_i$, either $v\geq {\rm top}(\mathfrak{w}_i)$ or ${\rm bot}(\mathfrak{w}_i)\geq v$, and $\mathfrak{w}_i\neq (v,v)$.
A standard sequence $\mathfrak{w}_1\geq\cdots\geq\mathfrak{w}_t$ of admissible pairs is called \textbf{anti-dominated} by $v$ if ${\rm bot}(\mathfrak{w}_t)\geq v$.
Let $SM^{v,v}$ denote the set of all $v$-compatible standard sequences of admissible pairs that are anti-dominated by $v$.
\edefn 
Example \ref{e.SMvv} below gives an illustration of $SM^{v,v}$.
\begin{ex}\label{e.SMvv}
{\rm
Let $d=4$ and $v=(1,2,3,5)$. Let $\mathfrak{w}_1=((2,4,6,8),(2,3,5,8))$ and $\mathfrak{w}_2=((1,2,4,6),(1,2,3,5))$. Then $\mathfrak{w}_1\geq\mathfrak{w}_2$ is an element of $SM^{v,v}$.}
\end{ex}
\bdefn\label{d.special-monomial}
A monomial $\mathfrak{S}$ of $mon{\mathfrak{N}(\gb)}$ is \textbf{special} if
\begin{enumerate}
    \item $\mathfrak{S}=\mathfrak{S}^{\hash}$ and
    \item the multiplicity of any diagonal element in $\mathfrak{S}$ is even.
\end{enumerate}
\edefn
Example \ref{eg.special} below gives an illustration of the special monomial.
\begin{ex}\label{eg.special}
{\rm
Let $d=4$ and $\gb=(1,2,5,6)$. Then the monomial $$\mathfrak{S}=\{(8,1),(8,1),(7,1),(8,2)\}$$ is a special monomial of $mon{\mathfrak{N}(\gb)}$.}
\end{ex}
Now we recall \cite[Proposition 4.1]{gr}, which has been used in the proof of the Theorem \ref{t.the-first-set}. \cite[Proposition 4.1]{gr} is stated below as Proposition~\ref{p.41ofgr}.
\begin{prop}\label{p.41ofgr}
There is a bijection between $SM^{\gb,\gb}$ and $mon{\mathfrak{N}(\gb)}$ that respects domination and degree. 
\end{prop}
Before we start the proof of the Theorem \ref{t.the-first-set}, let us recall the notation of $P$ and $\init S$, which will be used in the proof of the Theorem \ref{t.the-first-set}.
\bdefn\label{d.initial}
If $S$ is any nonempty subset of the polynomial ring $P:=K[X_{(r,c)}\st (r,c)\in\orootsb]$, such that $S \neq \{0\}$. We define $\init S$ to be the ideal $\langle \init(s)|s\in S\rangle$, where $\init$ is as in Definition~\ref{d.term-order}.
\edefn
\begin{thm}\label{t.the-first-set}
There exists a degree doubling injection from the set of all monomials of $P\setminus\init \mathcal{G}_{\ga,\gb}^{\gc}$ to the set of all non-vanishing special multisets on $\BBZ$ (bounded by $T_{\ga}$, $W_{\gc}$).
\end{thm}
\begin{proof}
Clearly, 
$$\init\mathcal{G}_{\ga,\gb}^{\gc}=\langle\init f_{\mathfrak{w},\gb}:\mathfrak{w}\ \text{is good}\rangle=\langle G^+\cup G^-\rangle,$$
where
$$G^+=\{x_{C^+}:C^+\ \text{a positive upper extended}\ \gb\text{-chain such that}\ C^+_{(1)}-C^+_{(2)}\nleq\gc\},$$
$$\mbox{and}\ G^-=\{x_{C^-}:C^-\ \text{a negative upper extended}\ \gb\text{-chain such that}\ \ga\nleq C^-_{(1)}-C^-_{(2)}\}.$$
Let 
$$G'^+:=\{x_{C^+}:C^+\ \text{a positive upper extended}\ \gb\text{-chain such that}\ C^+\nleq W_{\gc}\},$$
$$\mbox{and}\ G'^-:=\{x_{C^-}:C^-\ \text{a negative upper extended}\ \gb\text{-chain such that}\ T_{\ga}\nleq C^-\}.$$
It is then easy to observe that $G^+=G'^+$ and $G^-=G'^-$. Therefore
$$\init\mathcal{G}_{\ga,\gb}^{\gc}=\langle G'^+\cup G'^-\rangle.$$
The definition of a generating set for an ideal will now imply that $x_U$ is a monomial in $\init \mathcal{G}_{\ga,\gb}^{\gc}$ if and only if $x_U$ is a multiple of some $x_{C^+}$ or some $x_{C^-}$, where $C^+$ is a positive upper extended $\gb$-chain such that $C^+\nleq W_{\gc}$ and $C^-$ is a negative upper extended $\gb$-chain such that $T_{\ga}\nleq C^-$. Therefore $$x_U\ \mbox{is a monomial in}\ P\setminus\init \mathcal{G}_{\ga,\gb}^{\gc}$$
$$\Leftrightarrow x_U\ \text{is not divisible by any}\ x_{C^+}\ (\text{where}\ C^+\ \text{is a positive upper extended}\ \gb\text{-chain such that}\ C^+\nleq W_{\gc})$$
$$\text{or by any}\ x_{C^-}\ (\text{where}\ C^-\ \text{is a negative upper extended}\ \gb\text{-chain such that}\ T_{\ga}\nleq C^-)$$
$$\Rightarrow U\ \text{contains no extended upper}\ \gb\text{-chains}\ C\ \text{such that}\ T_{\ga}\nleq C^-\ \text{or}\ C^+\nleq W_{\gc}$$
Observe now that as the bijection of \cite[Proposition 4.1]{gr} respects domination, and \cite[Corollary 2.3.2]{papi1} holds true, so $C^+\nleq W_{\gc}$ implies that ${\rm top}(C^+)\nleq\gc$. A similar argument will show that $T_{\ga}\nleq C^-$ implies $\ga\nleq {\rm bot}(C^-)$. So we now have:
$$U\ \text{contains no extended upper}\ \gb\text{-chains}\ C\ \text{such that}\ T_{\ga}\nleq C^-\ \text{or}\ C^+\nleq W_{\gc}$$
$$\Rightarrow U\ \text{contains no extended upper}\ \gb\text{-chains}\ C\ \text{such that}\ \ga\nleq {\rm bot}(C^-)\ \text{or}\ {\rm top}(C^+)\nleq\gc$$
$$\Leftrightarrow \ga\leq {\rm bot}(C^-)\ \text{and}\ {\rm top}(C^+)\leq\gc\ \text{for any extended upper}\ \gb\text{-chain}\ C\ \text{in}\ U$$
$$\Leftrightarrow C\cup C^{\hash}\ \text{is bounded by}\ T_{\ga},W_{\gc}\ \text{for any extended upper}\ \gb\text{-chain}\ C\ \text{in}\ U,$$
where the last $\Leftrightarrow$ follows because ${\rm bot}(C^-)$ and ${\rm top}(C^+)$ are the two elements of $I(d)$ (as mentioned in Definition \ref{d.good-chains}) obtained by applying the map $BRSK$ to the monomial $C\cup C^{\hash}$, and the map $BRSK$ preserves domination.\par
Observe now that given any extended $\gb$-chain $D$ in $U\cup U^{\hash}$, we can naturally get hold of an extended upper $\gb$-chain $C$ (in $U$) from it in the following way:\\

If $D=(r_1,c_1)>\cdots>(r_t,c_t)$ and $(r_{i_1},c_{i_1}),\ldots,(r_{i_k},c_{i_k})$ (where $i_1<i_2<\cdots<i_k$) are such that $r_{i_j}>c_{i_j}^*$ for all $1\leq j\leq k$, then it is easy to check that the monomial formed by replacing all $(r_{i_j},c_{i_j})\ (1\leq j\leq k)$ in $D$ by $(c_{i_j}^*,r_{i_j}^*)$ forms an extended upper $\gb$-chain in $U$. Call this extended upper $\gb$-chain in $U$ as $C$.\par
Note that $D$ is an extended $\gb$-chain in the monomial $C\cup C^{\hash}$. So if the monomial $C\cup C^{\hash}$ is bounded by $T_{\ga},W_{\gc}$, then $T_{\ga}\leq D\leq W_{\gc}$.\par
Therefore\\
$$C\cup C^{\hash}\ \text{is bounded by}\ T_{\ga},W_{\gc}\ \text{for any extended upper}\ \gb\text{-chain}\ C\ \text{in}\ U$$
$$\Rightarrow T_{\ga}\leq D\leq W_{\gc}\ \text{for any extended}\ \gb\text{-chain}\ D\ \text{in}\ U\cup U^{\hash}$$
$$\Leftrightarrow U\cup U^{\hash}\ \text{is bounded by}\ T_{\ga},W_{\gc}.$$
The map $U\mapsto U\cup U^{\hash}$ from the set of all monomials of $P\setminus\init \mathcal{G}_{\ga,\gb}^{\gc}$ to the set of all non-vanishing special multisets on $\BBZ$ (bounded by $T_{\ga}$, $W_{\gc}$) is the required degree-doubling injection.
\end{proof}
Example \ref{e.the-first-theorem} below gives an illustration of Theorem \ref{t.the-first-set}.
\begin{ex}\label{e.the-first-theorem}
{\rm
Let $d=4$ and $\ga,\gb,\gc$ be as given in Example \ref{e.grobner-basis}. Here $T_{\ga}=\{(3,6)\}$ and $W_{\gc}=\{(3,1),(8,6)\}$. Take the monomial $U=\{(3,1),(3,2),(3,5)\}$ in $P\setminus\init \mathcal{G}_{\ga,\gb}^{\gc}$. It is now easy to verify that 
$$U\cup U^{\hash}=\{(3,1),(3,2),(3,5),(8,6),(7,6),(4,6)\}$$
is a non-vanishing special multiset on $\BBZ$ (bounded by $T_{\ga}$, $W_{\gc}$).}
\end{ex}
The following result follows easily from \cite[Propositions 6 and 7]{Br-La} followed by a proof similar to the proof of \cite[Proposition 3.9]{gr}:
\begin{prop}
The standard monomials on $Y_{\ga}^{\gc}(\gb)$ form a basis for $K[Y_{\ga}^{\gc}(\gb)]$. 
\end{prop}
\begin{thm}\label{t.the-second-set}
There exists a degree-halving injection from the set of all non-vanishing semistandard notched bitableaux on $(\BBZ)^*$ (bounded by $T_{\ga},W_{\gc}$) to the set of all standard monomials on $Y_{\ga}^{\gc}(\gb)$. 
\end{thm}
\begin{proof}
Given any non-vanishing semistandard notched bitableau $(P,Q)$ on $(\BBZ)^*$, let $P_1,\ldots,P_r$ (resp. $Q_1,\ldots,Q_r$) denote the rows of $P$ (resp. $Q$) from top to bottom. If $r$ is even (say, $r=2s$), then let us denote by 
$$\mu_1\leq\cdots\leq\mu_{2s}$$
the sequence $P_1-Q_1\leq\cdots\leq P_{2s}-Q_{2s}$. If $r$ is odd (say, $r=2s-1$), then let us denote by
$$\mu_1\leq\cdots\leq\mu_{2s}$$
the sequence $P_1-Q_1\leq\cdots\leq P_n-Q_n\leq\gb\leq P_{n+1}-Q_{n+1}\leq\cdots\leq P_{2s-1}-Q_{2s-1}$, where $(P_i,Q_i)$ (for $1\leq i\leq n$) is the negative part of $(P,Q)$, and $(P_j,Q_j)$ (for $n+1\leq j\leq 2m-1$) is the positive part of $(P,Q)$. We can then form the monomial $$f=f_{(\mu_2,\mu_1),\gb}f_{(\mu_4,\mu_3),\gb}\cdots f_{(\mu_{2s},\mu_{2s-1}),\gb},  $$ which belongs to $K[X_{(r,c)}|(r,c)\in\orootsb]$.\par 
The notched bitableau $(P,Q)$ is non-vanishing which implies that, any $\mu_i\ (1\leq i\leq 2s)$ which is not equal to $\gb$ is such that either $\mu_i<\gb$ or $\mu_i>\gb$. And this in turn implies that (\ref{e.standard2}), (\ref{e.standard3}), (\ref{e.standard4}) are satisfied by $f_{(\mu_2,\mu_1),\gb}\cdots f_{(\mu_{2s},\mu_{2s-1}),\gb}$. The notched bitableau $(P,Q)$ is semistandard $\Rightarrow\mu_1\leq\mu_2\leq\cdots\leq\mu_{2s-1}\leq\mu_{2s}$. Further, $(P,Q)$ is a notched bitableau on $(\BBZ)^*$ implies that the pairs $(\mu_2,\mu_1),\ldots,(\mu_{2s},\mu_{2s-1})$ are all admissible pairs. These two facts together imply that (\ref{e.standard1}) is satisfied by $f_{(\mu_2,\mu_1),\gb}\cdots f_{(\mu_{2s},\mu_{2s-1}),\gb}$. If in addition, $(P,Q)$ is bounded by $T_{\ga},W_{\gc}$, then it is implied that (\ref{e.standard5}) is satisfied by $f_{(\mu_2,\mu_1),\gb}\cdots f_{(\mu_{2s},\mu_{2s-1}),\gb}$. It is now easy to verify that $(P,Q)$ is a non-vanishing, semistandard, notched bitableau on $(\BBZ)^*$ (bounded by $T_{\ga},W_{\gc})\Rightarrow f_{(\mu_2,\mu_1),\gb}f_{(\mu_4,\mu_3),\gb}\cdots f_{(\mu_{2s},\mu_{2s-1}),\gb}$ is standard on $Y_{\ga}^{\gc}(\gb)$. Moreover, the degree of $(P,Q)$ equals the total number of boxes in $P$ (or $Q$). The total number of boxes in $P$ clearly equals $\Sigma_{i=1}^{r}|(P_i-Q_i)\setminus\gb|$, which in turn equals $\Sigma_{l=1}^{2s}|\mu_l\setminus\gb|$, which in turn equals twice the degree of $f_{(\mu_2,\mu_1),\gb}\cdots f_{(\mu_{2s},\mu_{2s-1}),\gb}$.\par
The map $(P,Q)\mapsto f_{(\mu_2,\mu_1),\gb}\cdots f_{(\mu_{2s},\mu_{2s-1}),\gb}$ is the required degree-halving injection.
\end{proof}
Example \ref{e.the-second-theorem} below illustrates Theorem \ref{t.the-second-set}.
\begin{ex}\label{e.the-second-theorem}
{\rm
Let $d=4$ and $\ga,\gb,\gc$ be as given in Example \ref{e.grobner-basis}. Here $T_{\ga}=\{(3,6)\}$ and $W_{\gc}=\{(3,1),(8,6)\}$. Let $$(P,Q)=\left(\young(3,4,37,38)\hspace{.2cm},\hspace{.2cm}\young(6,5,26,16)\right).$$
Now $(T_\alpha)_{(1)} - (T_\alpha)_{(2)} = (1,2,3,5)$, $(W_\gamma)_{(1)} - (W_\gamma)_{(2)} = (2,3,5,8)$, $P_1 - Q_1 = (1,2,3,5)$, $P_2 - Q_2 = (1,2,4,6)$, $P_3 - Q_3 = (1,3,5,7)$, and $P_4 - Q_4 = (2,3,5,8)$. \\
Since $P_i \ltS Q_i$ for $i = 1,2$ and $P_i \gtS Q_i$ for $i = 3,4$, so $(P,Q)$ is non-vanishing. Again $P_1 - Q_1 \leq P_2 - Q_2 \leq P_3- Q_3 \leq P_4-Q_4$, so $(P,Q)$ is semistandard. Also $(T_\alpha)_{(1)} - (T_\alpha)_{(2)} \ \leq P_1-Q_1$ and $P_4- Q_4 \ \leq (W_\gamma)_{(1)} - (W_\gamma)_{(2)}$, so $(P,Q)$ is bounded by $T_\alpha$ and $W_\gamma$.\par
Observe now that $(P,Q)$ is on $(\bar{\beta} \times \beta)^\star$ because:
\begin{enumerate}
\item $P_i = Q_i^\star$ for every row number $i$.
\item $(P,Q)$ does not contain any empty rows.
\item The total number of rows in $P$ (or $Q$) is 4, which is even.
\item The $\epsilon$-degrees of both $P_1-Q_1$ and $P_2-Q_2$ is 1. Also the $\epsilon$-degrees of both $P_3-Q_3$ and $P_4-Q_4$ is 2. That is, the $\epsilon$-degrees of $P_i-Q_i$ and $P_{i+1}-Q_{i+1}$ are equal for every $i$ odd.
\item The total number of box in P (or Q) is $6$ that is even.
\end{enumerate}
So $(P,Q)$ is on $(\bar{\beta} \times \beta)^\star$. Hence $(P,Q)$ is a non-vanishing semistandard notched bitableau on $(\BBZ)^*$ (bounded by $T_{\ga},W_{\gc}$). \\
Here $\mu_1=(1,2,3,5),\mu_2=(1,2,4,6),\mu_3=(1,3,5,7),\mu_4=(2,3,5,8)$. Under the degree-halving injective map of Theorem \ref{t.the-second-set}, $(P,Q)$ maps to $f_{(\mu_2,\mu_1),\gb}f_{(\mu_4,\mu_3),\gb}$, which is a standard monomial on $Y_{\ga}^{\gc}(\gb)$.}
\end{ex}

Before we start with the last theorem of this section, we state \cite[Lemma 6.3]{Kr-bkrs}, which is going to be used in the proof of the theorem. \cite[Lemma 6.3]{Kr-bkrs} is stated below as Lemma~\ref{l.lemma63ofkreiman}.
\begin{lem}\label{l.lemma63ofkreiman}
The map $BRSK$ is a degree-preserving bijection from the set of negative multisets on $\mathbb{N}^2$ to the set of negative semistandard notched bitableaux. 
\end{lem}
\begin{thm}\label{t.brsk-bijection}
The map $BRSK$ of \cite{Kr-bkrs} is a degree-preserving bijection from the set of all non-vanishing special multisets on $\BBZ$ (bounded by $T_{\ga},W_{\gc}$) to the set of all non-vanishing semistandard notched bitableaux on $(\BBZ)^*$ (bounded by $T_{\ga},W_{\gc}$).
\end{thm}
\begin{proof}
The fact that the map $BRSK$ of \cite{Kr-bkrs} is degree-preserving is obvious from \cite[Lemma 6.3]{Kr-bkrs} itself. For the rest of this proof, we will follow the notation and terminology of \cite[\S 4.1]{gr} as well as the notation and terminology of \cite{Kr-bkrs}.\par
There exists a natural injection from $SM^{\gb,\gb}$ to $\widetilde{SM^{\gb,\gb}}$ as given in \cite[\S 4.1]{gr}. Let $\widetilde{A^{\gb,\gb}}$ denote the image of $SM^{\gb,\gb}$ in $\widetilde{SM^{\gb,\gb}}$ under this injection. Let $\mathfrak{E}$ denote the set of all special monomials in $mon{\mathfrak{N}(\gb)}$. The map $\tilde{\pi}$ of \cite{Ko-Ra} is a degree and domination preserving bijection between the sets $\mathfrak{E}$ and $\widetilde{A^{\gb,\gb}}$. The set $\mathfrak{E}$ is the same as the set of all positive special multisets on $\BBZ$, and the set $\widetilde{A^{\gb,\gb}}$ is in a natural bijection (induced by the bijective map from $I_{\gb}^*$ to $I(d)$) with the set of all positive semistandard notched bitableaux on $(\BBZ)^*$. Also, we know from \cite[Corollary 2.3.2]{papi1} that the map $BRSK$ of \cite{Kr-bkrs} and the map $\tilde{\pi}$ of \cite{Ko-Ra} are the same on positive multisets on $\BBZ$. Moreover, it follows from \cite{Kr-bkrs} (see \cite[Lemma 7.2]{Kr-bkrs} and the fact that inside the proof of \cite[Proposition 2.3]{Kr-bkrs}, the inequality $b$ is actually an equality) that a positive multiset $U$ on $\BBZ$ is bounded by $\emptyset$, $W_{\gc}$ if and only if $BRSK(U)$ is bounded by $\emptyset$, $W_{\gc}$. Therefore, we can now conclude that the map $BRSK$ of \cite{Kr-bkrs} is a degree-preserving bijection from the set of all positive special multisets on $\BBZ$ (bounded by $\emptyset$, $W_{\gc}$) to the set of all positive semistandard notched bitableaux on $(\BBZ)^*$ (bounded by $\emptyset$, $W_{\gc}$).\par
The proof for the negative part is similar. For the negative part, the multisets (as well as the notched bitableaux) will be bounded by $T_{\ga}$, $\emptyset$ instead.
\end{proof}
Example \ref{e.brsk-bijection} below illustrates Theorem \ref{t.brsk-bijection}.
\begin{ex}\label{e.brsk-bijection}
{\rm
Let $d=4$ and $\ga,\gb,\gc$ be as given in Example \ref{e.grobner-basis}. Here $T_{\ga}=\{(3,6)\}$ and $W_{\gc}=\{(3,1),(8,6)\}$. Let 
$$\mathfrak{S}=\{(3,1),(3,2),(3,5),(8,6),(7,6),(4,6)\}.$$
Then $\mathfrak{S}$ is a non-vanishing special multiset on $\BBZ$ (bounded by $T_{\ga}$, $W_{\gc}$). 
$$BRSK(\mathfrak{S})=(P,Q)=\left(\young(3,4,37,38)\hspace{.2cm},\hspace{.2cm}\young(6,5,26,16)\right).$$
By Example \ref{e.the-second-theorem}, $BRSK(\mathfrak{S})=(P,Q)$ is a non-vanishing semistandard notched bitableau on $(\BBZ)^*$ (bounded by $T_{\ga},W_{\gc}$).}
\end{ex} 
\textit{The proof of Theorem \ref{t.main}} now follows easily from Theorems \ref{t.the-first-set}, \ref{t.the-second-set}, and \ref{t.brsk-bijection}.

\end{document}